# Combinatorial geometry takes the lead

**Andrei Okounkov**

**Abstract**

While the author is a professional mathematician, he is by no means an expert in the subject area of these notes. The goal of these notes is to share the author's personal excitement about some results of June Huh with mathematics enthusiasts of all ages, using maximally accessible, yet precise mathematical language. No attempt has been made to present an overview of the current state field, its history, or to place this narrative in any kind of broader scientific or social context. See the references in Section 9 for both professional surveys and popular science accounts that will certainly give the reader a broader and deeper understanding of the material.

## 1. Points, lines, and planes

Points and lines are the simplest geometric shapes and really primordial mathematical objects. Euclid opens his *Elements* by giving a definition of a point and a line, and his first postulate is that one can draw a straight line from any point to any point. While the need and standards for precise definitions in mathematics have only grown in the past $2.3 \cdot 10^3$ years, we imagine the reader has a good enough informal or formal grasp on lines and points to skip the definitions and focus on the basic geometric fact that two distinct point $P_1$ and $P_2$ determine a unique line through them. Somewhat unconventionally, we will denote this line $P_1 \vee P_2$.

While two points always lie on a line, three points $P_1, P_2, P_3$ may or may not be on a line. As we move the points around, the point $P_3$ is *typically* or *generically* not on the line $P_1 \vee P_2$, but in *special* cases it may be. The italicized words are important mathematical notions; we hope their meaning is intuitively clear.

Suppose we have $n = 3, 4, \ldots$ distinct points $P_1, \ldots, P_n$ in the plane, not all of them on same line. Generically, no three of these points will be on the same line, meaning that all lines $P_i \vee P_j$ will be distinct. Their number is thus the number of unordered pairs of numbers from $\{1, \ldots, n\}$, which can be computed as follows:

$$\# \text{ lines } = \frac{n(n-1)}{2} = 3, 6, 10, 15, \ldots, \quad n = 3, 4, 5, 6, \ldots .$$

See Figure (1) for illustration for $n = 7$. In particular, $n \ge 3$ generic points in the plane always determine $n$ or more lines.

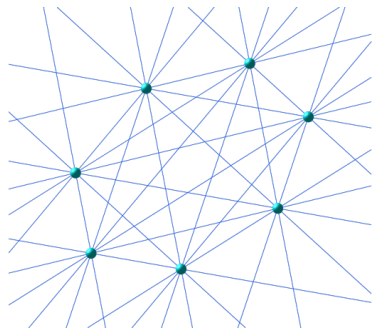

(1)

Let's see how the number of lines changes if we move the points into a special positition. For example, let's put $n - 1$ of them on a line, as in Figure (2). In this case, we get $n$ lines, so again at least as many lines as points.

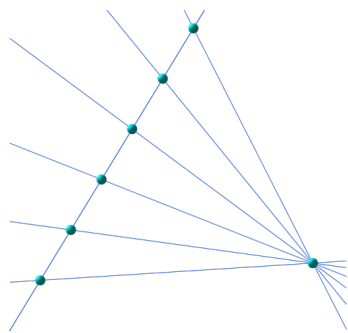

(2)

In general, it is classic result of de Brujin and Erdös [12] from 1948 that the number of lines determined by *n* points in the plane is at least *n*, unless all points lie on a single line.

Book 11 of the Elements opens with the definition of a solid, and Euclid proceeds with the development of the 3-dimensional geometry. Instead of *the* plane which previously contained the points $P_1, \ldots, P_n$, in three dimensions there are many planes and any triple of points $P_i, P_j, P_k$, not contained in a line, determines a unique plane $P_i \vee P_j \vee P_k$ that meets them.

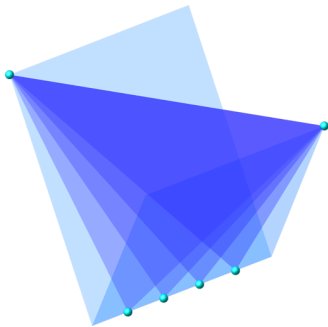

(3)

It is possible for *n* points to determine exactly *n* planes; see Figure (3) in which all but two points lie on a line. Theodore Motzkin [36] showed in 1951 that this indeed is the minimal possible number of planes[1].

## 2. Points, lines, planes, et cetera

It took a long time since Euclid for mathematicians and scientists to realize that it is both natural and important to study *d*-dimensional geometry for general $d = 1, 2, 3, 4, 5, \ldots$. A simple clear mathematical language of coordinates

$$\mathbb{R}^d = \{d\text{-tuples } (x_1, \ldots, x_d) \text{ of real numbers}\} \tag{4}$$

to describe the real *d*-dimensional space $\mathbb{R}^d$ was introduced in the 17th century by Fermat and Descartes. Tuples of numbers, sometimes with very large *d*, are abundant in both theoretical and applied contexts. But it was not until much later, less than 200 years ago, that the necessity and advantages of thinking about such *d*-tuples geometrically was realized.

A plane in a 3-dimensional space $\mathbb{R}^3$ is described by a linear equation

$$a_0 + a_1 x_1 + a_2 x_2 + a_3 x_3 = 0 \tag{5}$$

in which at least one of the coefficients $a_1$, $a_2$, or $a_3$ is not zero. Two sets of coefficients $(a_0, a_1, a_2, a_3)$ and $(a'_0, a'_1, a'_2, a'_3)$ determine the same plane if and only if

$$\begin{aligned}(a'_0, a'_1, a'_2, a'_3) &= c(a_0, a_1, a_2, a_3) \\ &= (ca_0, ca_1, ca_2, ca_3),\end{aligned} \tag{6}$$

for some nonzero number $c \neq 0$. Of course, multiplying an equation by a nonzero number does not change its solutions.

---

**1** In fact, Motzkin first conjectured this in his 1936 PhD thesis [35]!

A line in $\mathbb{R}^3$ is an intersection of two planes, thus the set of solution of a system of 2 linear equations. There are infinitely many planes containing a given line, and we can pick any two among them. In terms of the equations, this means that many transformations of a system of equations preserve their solutions. For instance, we can add to one of the equations any multiple of another equation.

Finally, a point $P$ in $\mathbb{R}^3$ is a solution of 3 linear equations, which we can choose to have the confusingly simple form

$$x_i = \text{the } i\text{th coordinate of } P, \quad i = 1, 2, 3. \tag{7}$$

In exactly the same fashion, a linear equation in $\mathbb{R}^d$ is said to determine a *hyperplane*, and points, lines, and *flats* of all other dimensions are described as the intersection of the corresponding number of hyperplanes, that is, as solutions of systems of linear equations. There is hardly anything more basic and fundamental in mathematics, science, technology, data analysis, et cetera, et cetera, than systems of linear equations. It is very likely that many or most readers of these notes have met them before. Those who would like a reminder or an explanation will find it in Appendix A.

The basic geometric facts like:

2 points $P_1, P_2$, when distinct, lie on a unique line $P_1 \vee P_2$,

3 points $P_1, P_2, P_3$, not contained in a line, lie in a unique plane $P_1 \vee P_2 \vee P_3$,

...

$r$ points $P_1, \ldots, P_r$, not contained in a $(r-2)$-dimensional flat,

$$\text{lie in a unique } (r-1)\text{-dimensional flat } P_1 \vee P_2 \vee \cdots \vee P_r, \tag{8}$$

continue to hold in any dimension $d$. The minimal flat containing some points $P_1, \ldots, P_k$ will be denoted $P_1 \vee P_2 \vee \cdots \vee P_k$ and called the *span* of these points.

It is natural to ask how many flats of each dimension can $n$ points in $\mathbb{R}^d$ determine. Since it takes $r$ points to determine an $(r-1)$-dimensional flat, we will define the *rank* of such flat to equal $r$.

For instance, $n$ generic points $P_1, P_2, \ldots, P_n$ determine

$$\binom{n}{2} \text{ lines}, \binom{n}{3} \text{ planes}, \ldots, \binom{n}{r} \text{ rank } r \text{ flats}, \ldots \tag{9}$$

because we can choose so many $r$-element subsets from an $n$-element sets. Here

$$\binom{n}{r} = \frac{n!}{r!(n-r)!}, \quad n! = 1 \cdot 2 \cdot 3 \cdot \cdots \cdot (n-1) \cdot n,$$

denote the binomial coefficients $\binom{n}{r}$ and the *factorial* $n!$ of $n$. Factorials and binomial coefficients are as fundamental to mathematics and as ancient as points and lines, appearing in very old Indian, Persian, and Chinese texts long before becoming known in Europe in late Renaissance.

Many beautiful elementary properties of the binomial coefficients inspire combinatorialists to look for similar patterns in other, more complicated, sequences of numbers. In ([9](#)) we get some initial segment of the binomial coefficients, which look as follows for $d = 6$, $n = 10$ and $r = 1, \ldots, d$

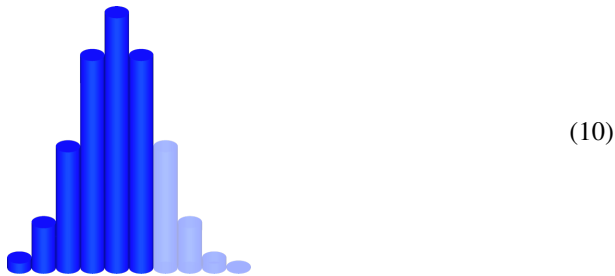

(10)

The transparent columns here represent the unused binomial coefficients with $r > d$. Two properties of this sequence of numbers are apparent. First, it is *unimodal*, that is, the numbers first increase and then decrease. Second, it is *top-heavy*, which can be quantified as

$$\binom{n}{r} \leq \binom{n}{d-r}, \quad \text{provided} \quad 2r \leq d \leq n. \tag{11}$$

So far, this was about $n$ generic points in $\mathbb{R}^d$. Now let us allow the points $P_1, P_2, \ldots, P_n$ to be in some special position (and there are a great many ways in which a point configuration can be special for $n$ large). Let $\mathscr{F}_r$ denote the set of rank $r$ flats determined by the $P_i$'s. In particular,

$$\mathscr{F}_1 = \{P_1, \ldots, P_n\}, \tag{12}$$

while $\mathscr{F}_2$ are the lines in $\mathbb{R}^d$ containing at least two of the $P_i$'s. Generalizing what we have seen for generic points, Rota [39] conjectured the *unimodality* of the sequence $|\mathscr{F}_r|$, where $|\mathscr{F}_r|$ denotes the numbers of elements, or *cardinality*, of $\mathscr{F}_r$. Dowling and Wilson [15, 16] conjectured that the sequence $|\mathscr{F}_r|$ is top-heavy. These questions remained open for a very long time, but now the top-heavy conjecture and the increasing part of the unimodality conjectures are proven as a corollary of a theorem of June Huh and Botong Wang that will be discussed in the next section.

Why is the top-heavy conjecture so interesting? "*It indicates a deep hidden reciprocity!*", says Gil Kalai who presented June Huh's Fields Medal laudatio at ICM 2022. June Huh says he became interested in the top-heavy conjecture as a result of being intrigued by the "top-heavy phenomena" for lower Bruhat intervals in Coxeter groups that are proved using Elias–Williamson's combinatorial Hodge theory for Soergel bimodules. A curious reader will find out what this is about in the references [9, 17, 24, 34].

### 3. Matching flats to flats

Suppose we want to prove that one set, such as $\mathscr{F}_r$, has fewer elements than some other set, such as $\mathscr{F}_{r'}$. These sets may be complicated and the exact counts of elements in each of them may be hard to perform. However, we may be able to prove the inequality between $|\mathscr{F}_r|$ and $|\mathscr{F}_{r'}|$ without actually doing either count. It suffices to assign to each element $F \in \mathscr{F}_r$ an element $\iota(F) \in \mathscr{F}_{r'}$ so that distinct $F_1 \neq F_2$ are assigned distinct $\iota(F_1) \neq \iota(F_2)$.

Mathematicians have special words for any procedure $\iota$ that assigns an element $\iota(F)$ of some "target" set like $\mathscr{F}_{r'}$ to an element $F$ of some "source" set like $\mathscr{F}_r$. We say that $\iota$ is a function or a map from $\mathscr{F}_r$ to $\mathscr{F}_{r'}$ and write:

$$\iota : \mathscr{F}_r \to \mathscr{F}_{r'} . \tag{13}$$

When (13) takes distinct elements to distinct elements, we say that $\iota$ is *injective* or one-to-one. An injective map between two finite sets exists if and only if the cardinality of the source is less than or equal to the cardinality of the target.

Conversely, a map is called *surjective* or onto, if every element in the target is assigned to some element of the source. A surjective map implies the opposite inequality between the cardinalities of the two sets. A schematic example of an injective and a surjective set from a set of circles to a set of stars can be seen in Figure (14).

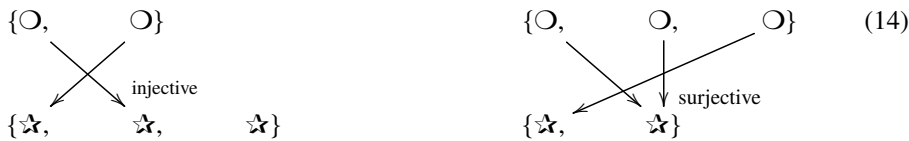
(14)

For an injective map, every star is the target of $\leq 1$ arrows; for a surjective map, every star is the target of $\geq 1$ arrows.

Since the source and the target in (13) have a geometric meaning, we can ask for the map $\iota$ to reflect this geometric meaning. It is nice to require that the flat $\iota(F)$ contains the flat $F$ for all $F$. We will call such assignment a *matching*. In Figure (15), the reader can see examples of a non-injective and an injective matching between points and lines from Figure (2).

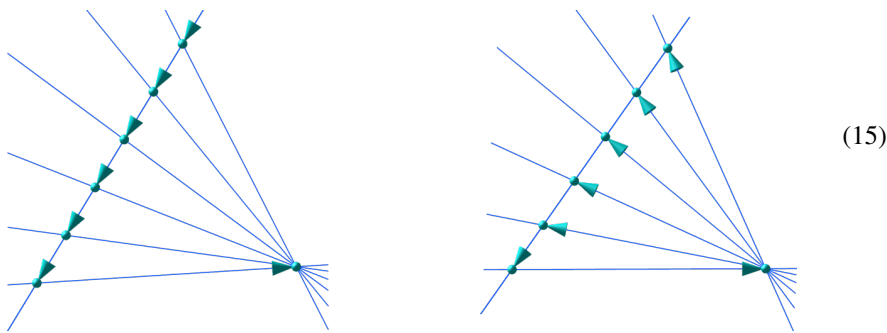
(15)

In Figure (15), we give each point a conic tail in the direction of the matched line.

These notes are about some very general and powerful results proved by June Huh and his collaborators Tom Braden, Jacob P. Matherne, Nicholas Proudfoot, and Botong Wang. We will be stating and explaining them in ascending generality, starting with the following most basic version:

**Theorem 1** ([22]). *For any n-tuple of points $P_1, \ldots, P_n \in \mathbb{R}^d$ not contained in a hyperplane, there exists an injective matching* (13) *from rank r flats to rank r' flats spanned by these points provided $r \leq r'$ and $r + r' \leq d + 1$.*

For instance, in the plane $d = 2$, the only interesting case is $r = 1$ and $r' = 2$. Theorem 1 then says that for every point one can chose a line containing it, in such a way that different points are assigned different lines. In $\mathbb{R}^3$ we can have $(r, r') = (1, 2)$ or $(r, r') = (1, 3)$, which means that for every point we can choose a line, and also a plane containing it. In every dimension, Greene showed that points can be matched to the hyperplanes they define [19].

It is not so clear at present what happens outside the $r + r' \leq d + 1$ range, including the decreasing part of the unimodality conjecture. For instance, Dilworth and Greene constructed in [14] a configuration of 21 points in a 10-dimensional space such that there is no injective or surjective matching $\mathscr{F}_6 \to \mathscr{F}_7$.

## 4. Rank and Matroids

Suppose that for some concrete collection $P_1, \ldots, P_n$ we want a computer program to either construct or verify an injective matching described in Theorem 1. Or maybe we would like to experiment in the range $r + r' > d + 1$. Whatever our goals, we will need the program to manipulate the information about the position of the points $P_1, \ldots, P_n$. It should be able to either determine or remember which subsets

$$S \subset \{P_1, \ldots, P_n\}$$

of points lie on a line, in a plane, et cetera. So it reasonable to think that in our program there should be a procedure that either computes or looks up the function

$$\text{rank}(S) = \dim \text{span}(S) + 1 \, . \tag{16}$$

It is easy to see that all other notions discussed so far can be easily expressed in terms of the function (16). For instance, a subset $S$ corresponds to a flat of rank $r$ if and only rank$(S) = r$ and

$$\text{rank}(S \cup P_i) = \text{rank}(S) + 1 \, , \quad \text{for any } P_i \notin S \, . \tag{17}$$

In other words, flats $F \in \mathscr{F}_r$ corresponds to subsets $S$ of rank $r$ that are maximal with respect to inclusion.

As we change the position of the points $P_1, \ldots, P_n$, the corresponding rank functions will also change, but they will always satisfy the equalities:

$$\text{rank}(\{P_i\}) = 1 \, , \quad i = 1, \ldots, n \, , \tag{18}$$

and the inequalities

$$\text{for any } S \subset S', \quad \text{rank}(S) \leq \text{rank}(S') \, , \tag{19}$$

$$\text{for any } S_1, S_2, \quad \text{rank}(S_1 \cup S_2) \leq \text{rank}(S_1) + \text{rank}(S_2) - \text{rank}(S_1 \cap S_2) \, . \tag{20}$$

The intersection $S_1 \cap S_2$ in (20) may be the empty set $\varnothing$, and by definition

$$\text{rank}(\varnothing) = 0 \, . \tag{21}$$

The geometrically obvious inequality rank$(S) \leq |S|$ is then a formal corollary of (20).

One may wonder whether, in fact, the above properties characterize all possible rank functions for points in a space of dimension $d$, where $d = \text{rank}(\{P_1, \ldots, P_n\}) - 1$? And maybe a proof of Theorem 1 may be found by exploring formal consequences of (19) and (20)? Turns out, the answers to these questions are emphatic "no" and "yes", respectively.

The above properties of the rank function give one of the many equivalent axiomatic definitions of a matroid[2]. Matroids were introduced by Hassler Whitney in 1935 as combinatorial generalization of incidence relations between flats of different dimensions, and have since found an abundance of applications across mathematics and computer science, both pure and applied. "*Matroid theory is a triumph in the pursuit of both abstraction and concrete simple examples*", says Gil Kalai.

While we will discuss a few examples below, it should be made very clear now that matroids constitute a very rich and diverse universe, much larger than what we will explore in these notes. This makes the following result of Tom Braden, June Huh, Jacob Matherne, Nicholas Proudfoot, and Botong Wang extremely remarkable and powerful

**Theorem 2.** [10] *The injective matching*

$$\iota : \mathscr{F}_r \to \mathscr{F}_{r'}, \quad r \leq r',$$

*as in Theorem 1, exists for flats of any matroid M provided $r + r' \leq \text{rank}(M)$.*

While the definition of a matroid is very short, the argument leading to the proof of Theorem 2 is very, very complex. The best we can hope to do in these introductory notes, is to explain some earlier results and ideas in various areas of mathematics that may be listed among the precursors and inspirations for the fantastic achievement of [10].

At several points in our narrative, we will be coming back the following extraordinary feature of Theorem 2. In geometry, there is a constant dialog between the continuous and the discrete. Of course, there is a fundamental unity in mathematics and good mathematics is constantly transcending apparent boundaries between different subfields. Still, there is a clear difference between a matroid, which is combinatorial abstraction of a geometric configuration, and objects like a geodesic on a manifold, a minimal surface, or a harmonic form that require noncombinatorial methods to define and study.

One can compare and contrast the continuous and the discrete in many different ways, but one fundamental difference is the presence of *limits* in the continuous world. Of course, limits are absolutely central to mathematics and many crucial mathematical objects, like the exponential function $e^x$, are transcendental in the sense that a limit is required to define or compute it. It is, however, an interesting question how much extra mileage one can get from using analytic tools to investigate combinatorial objects. As we will see, at the heart of Theorem 2, lies a certain *hard Lefschetz* property, which for many years was firmly associated with continuous, noncombinatorial geometry.

---

[2] More precisely, the condition (18) means that here we focus on so-called loopless matroids. Given a rank function, one defines the flats of a matroid as in (17). Conversely, the rank function may be reconstructed from the data of the flats.

8        Andrei Okounkov

Theorem 2 applies to an *arbitrary* matroid, a purely combinatorial object. This is already very remarkable. But what is really amazing is the combinatorial and algebraic framework built in [10] to prove Theorem 2. It produces the required hard Lefschetz property from purely combinatorial, finite ingredients.

We will come back to these points later in the narrative. First, in the next section we want to discuss some examples of matroids beyond what we have seen so far. We warn the reader that these examples still cover a vanishing fraction of the universe of matroids.

## 5. Some examples of matroids
### 5.1. Points in $\mathbb{F}^d$, where $\mathbb{F}$ is a field

The first generalization concerns the coordinates in (4). There we took a $d$-tuple of real numbers, while instead we could have taken $x_i$ to be rational numbers $x_i \in \mathbb{Q}$, or complex numbers $x_i \in \mathbb{C}$, or elements in an arbitrary field.

In mathematics, a field $\mathbb{F}$ is a set with special elements $0, 1 \in \mathbb{F}$ and binary operations $+, -, \times, /$ obeying all the usual laws of arithmetic for rational numbers $\mathbb{Q}$. An attentive reader will notice that the division is not really binary, since one cannot divide by zero. Instead, it is a function of the form

$$\mathbb{F} \times (\mathbb{F} \setminus \{0\}) \xrightarrow{(a,b) \mapsto a/b} \mathbb{F}\,.$$

For a dramatic example, we can take $\mathbb{F}_2 = \{0, 1\}$. Since $a + 0 = a$, $a \times 0 = 0$, and $a \times 1 = a$ in any field, most sums and products are already specified. The only interesting one is $1 + 1 = ?$ We invite the reader to check that

$$1 + 1 = 0 \qquad (22)$$

is the only logical option and this indeed defines a field with two elements.

All constructions of Section 2 extend verbatim to any field. Note, however, that some configurations of points and lines can be realized using $\mathbb{F}_2$ and cannot be realized with real numbers, see (33) for an example.

Inspired by (22), we can ask when it is possible that

$$\underbrace{1 + 1 + \cdots + 1}_{p \text{ times}} = 0 \qquad (23)$$

in a field $\mathbb{F}$. Minimal possible $p$ with this property should be a *prime*, as we invite the reader to check. It is called the *characteristic* of the field $\mathbb{F}$. An example of a field of characteristic $p$ is given by the residues

$$\mathbb{F}_p = \{0, 1, \ldots, p-1\}\,.$$

modulo a prime $p$. The operations are defined by

$$(a, b) \xrightarrow{+} (a+b) \bmod p\,, \quad (a, b) \xrightarrow{\times} ab \bmod p\,, \qquad (24)$$

where $a + b$ and $ab$ are the usual addition and multiplication of integers. But since there is no division of integers, the existence of multiplicative inverse requires a minute of thought. For instance, one may notice that the multiplication by $a$ map

$$\mathbb{F}_p \ni b \mapsto ab \in \mathbb{F}_p \qquad (25)$$

is injective for $a \ne 0$, and therefore also surjective, as the source and target sets have the same cardinality. A map that is both injective and surjective is called *bijective*; those are the maps that have an inverse map. The inverse to multiplication by $a$ is, by definition, the division by $a$. For fun, the reader may want to compute the inverses of elements in $\mathbb{F}_p$ for $p = 3, 5, 7$.

Mathematicians usually think about fields in terms of field *extensions*. Concretely, one describes new fields $\mathbb{F}'$ in terms of some previously understood subfield $\mathbb{F} \subset \mathbb{F}'$ and the new elements $x_1, x_2, \dots \in \mathbb{F}'$ that have to be added to $\mathbb{F}$ to generate all elements of $\mathbb{F}'$ by arithmetic operations. One writes $\mathbb{F}' = \mathbb{F}(x_1, x_2, \dots)$ to denote this situation. For example

$$\mathbb{C} = \mathbb{R}(\sqrt{-1}) \, .$$

All information about the field extension is contained in the polynomial equations satisfied by the elements $x_1, x_2, \dots$ with coefficients in $\mathbb{F}$. For example, the element $i = \sqrt{-1}$ satisfies the equation

$$i^2 + 1 = 0 \, . \qquad (26)$$

Using this equation, we simplify powers $i^k$ and, in particular, compute the product of two complex numbers as follows

$$(a_1 + a_2 \cdot i)(b_1 + b_2 \cdot i) = (a_1 b_1 - a_2 b_2) + (a_1 b_2 + a_2 b_1) \cdot i \, . \qquad (27)$$

In parallel to (25), we invite the reader to check that this is injective and surjective for $a_1^2 + a_2^2 \ne 0$ and compute $(a_1 + a_2 \cdot i)^{-1}$. After this exercise, the reader may want to find the formula for the inverse in $\mathbb{Q}(\sqrt{2})$.

To reiterate, while field theory is ultimately about solutions of polynomial equations, it is much more effective to *use* equations to learn about their solutions as opposed to "solving" the equations in the sense of looking for some complicated formulas giving the solutions in terms of the coefficients. We took this little detour into algebra now because later it will be crucial to *use* certain algebraic equations to deduce the information about flats in a matroid.

### 5.2. Projective spaces

While two points in a plane $\mathbb{R}^2$ always determine a line, two lines in a plane usually intersect in a single point, but not always. Sometimes lines can become parallel and then their point of intersection runs off to infinity. Projective geometry adds these points at infinity to the plane $\mathbb{R}^2$ so that two lines always intersect at a point. Many other geometric statements no longer require considering various cases, either[3].

---

    **3**    For instance, the hyperbola, the parabola, and the ellipse are the same geometric shapes in projective geometry !

**10**    Andrei Okounkov

The $d$-dimensional projective space over an arbitrary field $\mathbb{F}$ is easy to define using coordinates as follows

$$\mathbb{P}^d(\mathbb{F}) = \{(d+1)\text{-tuples } (x_0, x_1, \ldots, x_d) \text{ of elements of } \mathbb{F},$$
$$\text{not all zero, up to proportionality } \sim\}, \tag{28}$$

where proportionality means

$$(x_0, x_1, \ldots, x_d) \sim c(x_0, x_1, \ldots, x_d), \quad c \in \mathbb{F} \setminus \{0\}. \tag{29}$$

Recall that we already met such identification of proportional tuples in (6) when we talked about hyperplanes in $\mathbb{R}^d$.

The $d$-dimensional space $\mathbb{F}^d$ is naturally embedded in the projective space as the set where $x_0 \neq 0$. Indeed, when $x_0 \neq 0$ we can choose a unique $c$ in (29) to make $x_0 = 1$, and so we get

$$\mathbb{F}^d = \{(1, x_1, \ldots, x_d)\} \subset \mathbb{P}^d(\mathbb{F}). \tag{30}$$

The points with $x_0 = 0$ are the points "at infinity". They form a smaller projective space $\mathbb{P}^{d-1}(\mathbb{F})$.

By definition, a hyperplane in $\mathbb{P}^d(\mathbb{F})$ is defined by an equation of the form

$$a_0 x_0 + a_1 x_1 + \cdots + a_d x_d = 0, \tag{31}$$

in which some of the coefficients $a_i \in \mathbb{F}$ are not zero. In particular, the "infinity" is the hyperplane $x_0 = 0$.

Note that (31) is unchanged if we scale all variables $x_i$ or all variables $a_i$ by some constant $c \neq 0$. Thus the hyperplanes in $\mathbb{P}^d(\mathbb{F})$, as described by their coefficients $(a_0, \ldots, a_d)$, form another projective space, called the *dual* projective space. This basic duality underlies many remarkable facts in geometry and combinatorics.

If the field $\mathbb{F}$ is finite, one can get very interesting matroids by taking *all* points of $\mathbb{P}^d(\mathbb{F})$ as the points $P_1, \ldots, P_n$. For example, if $\mathbb{F} = \mathbb{F}_2$, there is no need to worry about proportionality (29), and so we get 7 points

$$\mathbb{P}^2(\mathbb{F}_2) = \{(1,0,0), (0,1,0), (0,0,1), (1,1,0), (0,1,1), (1,0,1), (1,1,1)\}. \tag{32}$$

By duality, there are 7 hyperplanes in $\mathbb{P}^2(\mathbb{F}_2)$. Each of them, in some coordinates, represents $\mathbb{P}^1(\mathbb{F}_2)$, and hence contains 3 points from (32). We invite the reader to check that the resulting configuration of points and lines looks as follows:

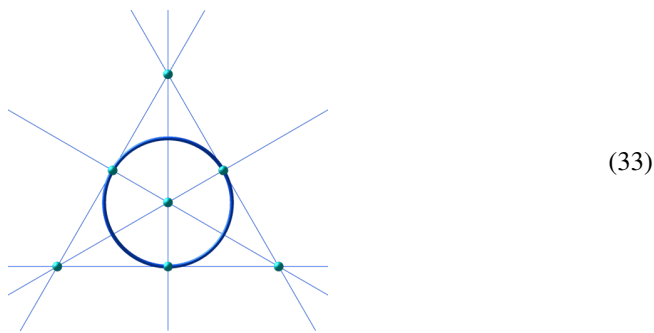

(33)

People who have read about the math behind the game of Dobble/Spot it! certainly recognize this picture. In (33), the line

$$x_1 + x_2 + x_3 = 0$$

is plotted as a circle precisely because this line does not meet the points (32) with real coefficients. With real coefficients, the three points on the circle define three different lines.

In fact, the matroid (33) can be realized in $\mathbb{P}^2(\mathbb{F})$ for some field $\mathbb{F}$ if and only if the characteristic of the field $\mathbb{F}$ equals 2, as noted already by Whitney. In general, whether a given matroid can be put into $\mathbb{P}^d(\mathbb{F})$, or as one says can be *linearly realized* over $\mathbb{F}$, is an interesting and important question.

Projective geometry is a very classical and very beautiful subject. Any result about incidence of points, lines, etc. in projective geometry is a potential obstruction to linear realizability of a given matroid over any field. For instance, the Pappus theorem says that if the three top points in (34) are collinear, and the three bottom points are also collinear, then so are the three middle points. The line which Pappus proved exists, is highlighted in (34)

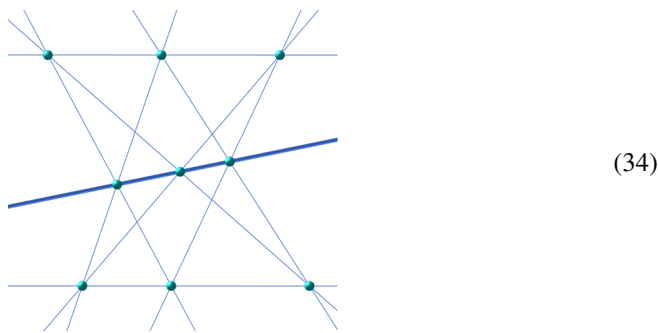

(34)

One can, however, declare the three middle points to be noncollinear without violating the matroid axioms. This gives a concrete example of a matroid that is not linearly realizable. In fact, the vast majority of matroids are such [37]. One way to think about them is to take a realizable matroid and add/remove some flats making sure the axioms are still obeyed, like we just did with with the non-Pappus matroid.

### 5.3. Field extensions

This class of examples is a little more advanced, so may be skipped on the first reading. Let $\mathbb{F}$ be a field and consider an extension

$$\mathbb{F}' = \mathbb{F}(y_1, \ldots, y_n).$$

Note that we have $n$ elements generating the extension, whereas earlier $n$ was the number of points. This is intentional and we will define a rank function on subsets $S \subset \{1, \ldots, n\}$ by

$$\text{rank } S = \text{transcendence degree of } \mathbb{F}(\{y_i\}_{i \in S}) \text{ over } \mathbb{F}. \tag{35}$$

Here the transcendence degree of a field over a subfield $\mathbb{F}$ is the maximal number $t$ of elements $f_1, f_2, \ldots, f_t$ that satisfy no polynomial equation $P(f_1, \ldots, f_t) = 0$ with coefficients in $\mathbb{F}$. The rank function (35) satisfies the matroid axioms.

Algebraic geometers associate a geometric image to this algebraic definition by thinking about

$$\mathbb{F}' = \text{rational functions on } Y, \quad Y \subset \mathbb{F}^n, \tag{36}$$

where $y_1, \ldots, y_n$ are the coordinate functions on $\mathbb{F}^n$ and $Y$ is an irreducible algebraic variety. By definition, an algebraic variety is defined by polynomial equations and it is irreducible if it is not a union of two other varieties. In (36) we identify two rational functions $f(y)$ and $f'(y)$ if they are equal on $Y$. With such reformulation

$$\text{rank } S = \dim \left( \text{projection of } Y \text{ to the coordinates } \{y_i\}_{i \in S} \right). \tag{37}$$

Matroids that can be put in this form are called *algebraically realizable*. For instance, the non-Pappus matroid on 9 points is algebraically realizable over $\mathbb{F} = \mathbb{F}_2$; see [32].

The reader should check that a linear realization yields an algebraic realization with $Y \subset \mathbb{F}^n$ being a linear subspace. In this case, it was shown in [3] that the matroid (35) controls many important properties of the closure

$$\overline{Y} \subset \left( \mathbb{P}^1(\mathbb{F}) \right)^n \tag{38}$$

of the linear space $Y$ in a product of projective lines over $\mathbb{F}$. We will come back to (38) in Section 8.3.

In another direction, Ingleton notes [25] that the tangent space to the generic points of $Y$ provides a linear representation of an algebraic matroid, provided the characteristic of $\mathbb{F}$ is zero. This does not work in positive characteristic as the non-Pappus matroid demonstrates[4].

Some simple explicit matroids can be shown to be algebraically nonrealizable; see [26].

### 5.4. Tropical realization of matroids

So far, we have looked at different classes of examples of matroids, always stressing the fact that these do not cover the great diversity of the world of matroids. Very remarkably, however, there is a class of examples that gives *all* matroids. This is the case for the *tropical* analog of the construction from Section 5.3, in which it is enough to take $Y$ to be a tropical linear space.

This was discovered by Bernd Sturmfels in [44]. See Appendix C for a few introductory comments, [27, 30, 33] for a proper introduction to the subject, and [1, 2, 6–8, 18] for a sample of exciting recent advances in this direction. June Huh says: "*Mathematicians discovered tropical varieties by tropicalizing algebraic varieties, but only a tiny fraction of tropical varieties are tropicalizations of algebraic varieties. Tropical varieties are geometric objects that try to teach us a new kind of geometric intuition through their diversity.*"

In the spirit of this narrative, one should wish the best of success to all present and future combinatorial geometers in extending classical geometric results to this combinatorial setting. It is both very beautiful and important for applications.

---

**4**  In characteristic $p$ one struggles with tangent spaces due to the fact that $(x^p)' = px^{p-1} = 0$.

## 6. Graded Möbius algebra

A certain algebraic language will be required the capture the essense of that is happening in Theorem 2. Most importantly, we will need explain one more meaning that the mathematicians attach to the word *algebra*.

### 6.1. Algebras

Consider a field extension $\mathbb{F}' = \mathbb{F}(x)$ generated by one element satisfying a polynomial equation of degree $d$ with coefficients in $\mathbb{F}$. For instance, it can be $\mathbb{F}' = \mathbb{Q}(\sqrt[4]{-2})$, which means that the coefficients $\mathbb{F} = \mathbb{Q}$ are rational numbers and the new element $x$ satisfies the equation

$$x^4 + 2 = 0. \tag{39}$$

We can think of $x^4 = -2$ as a substitution rule that we can instruct our computer to apply any time it sees a power $x^k$ with $k > 3$. Using this substitution rule, we can describe $\mathbb{F}'$ by 4-tuples of rational numbers

$$\mathbb{F}' = \{a_0 + a_1 x + a_2 x^2 + a_3 x^3, a_i \in \mathbb{Q}\}, \tag{40}$$

and thus we can picture $\mathbb{F}'$ as a 4-dimensional linear space[5] $\mathbb{Q}^4$. To multiply two general elements of $\mathbb{F}'$

$$(a_0 + a_1 x + a_2 x^2 + a_3 x^3)(b_0 + b_1 x + b_2 x^2 + b_3 x^3) = ?,$$

we have to expand out and use the substitution rule (39). This rule can phrased as follows

$$(a_0 + a_1 x + a_2 x^2 + a_3 x^3)(b_0 + b_1 x + b_2 x^2 + b_3 x^3) =$$
$$= \underbrace{(c_0 + c_1 x + c_2 x^2 + c_3 x^3)}_{\text{product in } \mathbb{F}'} + \underbrace{\text{something} \cdot (x^4 + 2)}_{\text{discard}}, \tag{41}$$

where something refers to some polynomial in $x$. To see that multiplication by a nonzero element is invertible, it suffices to check $\mathbb{F}'$ has no nontrivial divisors of zero. (We already used this logic when we were inverting (25).) But any $a_0 + a_1 x + a_2 x^2 + a_3 x^3$ that divides zero in $\mathbb{F}'$ will have to divide $x^4 + 2$, whereas this polynomial cannot be nontrivially factored into polynomials with rational coefficients[6]. (Check this!)

What if we replaced 2 by 0 in (39), that is, what if we used a simpler equation $x^4 = 0$? The presentation (40) would still be valid and the multiplication would take a simpler form

$$(a_0 + a_1 x + a_2 x^2 + a_3 x^3)(b_0 + b_1 x + b_2 x^2 + b_3 x^3) =$$
$$= \underbrace{(c_0 + c_1 x + c_2 x^2 + c_3 x^3)}_{\text{product}} + \underbrace{\text{terms of degree} \geq 4 \text{ in } x}_{\text{discard}}, \tag{42}$$

---

[5] Readers who would like a bit more details about linear spaces will find them in Section A.3.
[6] Mathematicians say (39) is *irreducible*. We have already used this term in exactly this meaning Section 5.3.

meaning that

$$c_0 = a_0 b_0,$$
$$c_1 = a_1 b_0 + a_0 b_1,$$
$$c_2 = a_2 b_0 + a_1 b_1 + a_0 b_2, \quad \text{et cetera}.$$

Of course, this will no longer be a field, because multiplication by $x$ is not invertible. But it is still a viable algebraic object that we will denote by a different letter

$$\begin{aligned} \mathbb{A} &= \{a_0 + a_1 x + a_2 x^2 + a_3 x^3, a_i \in \mathbb{Q}\} \\ &= \mathbb{Q}[x]/(x^4 = 0), \end{aligned} \tag{43}$$

lest somebody thinks it is a field.

We see that $\mathbb{A}$ is a linear space over a field $\mathbb{F}$ that has a product operation satisfying all the rules of arithmetic except those involving division. Mathematicians call such an object an *algebra*[7], not to be confused with algebra as an area of mathematics that studies fields, algebras, and many other important structures. To distinguish between field and algebras, mathematicians use square brackets in (43) and they also like to to write the equations imposed on $x$ as in (43). If there were no equations on $x$, the we would simply get the algebra of polynomials in $x$ with coefficients in $\mathbb{Q}$. That algebra is denoted $\mathbb{Q}[x]$.

The reader may wonder what could be the purpose of studying equations like $x^4 = 0$. Doesn't this just mean that $x = 0$? In fact, no, and there are very natural geometric situation where relations like this appear. Let's look at the following table

$$\begin{array}{c|c|c|c|c} 1 & x & x^2 & x^3 & x^4 = 0 \\ \hline \text{space} & \text{plane} & \text{line} & \text{point} & \varnothing \end{array}, \tag{44}$$

and note that the following parallels

$x \cdot x = x^2$      two general planes in space intersect in a line,

$x \cdot x^2 = x^3$      a plane and a general line in space intersect in a point,

$x \cdot x^3 = 0$      a plane and a general point in space intersect in an empty set.

Note that the 3-space here can be over an arbitrary field, which has nothing to do with the rational coefficients we had in (43). For a $d$-dimensional space, we should replace $x^4 = 0$ by $x^{d+1} = 0$. We will come back to these parallels in Section 8, but for now notice the potential for algebras to encode combinatorial information.

### 6.2. Graded algebras

This potential to encode combinatorial information gets amplified when we consider graded algebras. Let

$$\mathbb{A} = \mathbb{F}[x_1, \ldots, x_N] \big/ (\text{relations})$$

---

[7]      Or a *commutative algebra* to be more precise, since the product in $\mathbb{A}$ still obeys the commutative law of the arithmetic.

15      Combinatorial geometry takes the lead be an algebra generated by generators $x_i$ subject to some relations. By definition, $\mathbb{A}$ is graded if every generator $x_i$ is assigned a positive integer $\deg x_i = 1, 2, \ldots$ so that all relations only involve monomials of the same total degree in $x_1, \ldots, x_N$. For instance, $x^4 = 0$ is a good relation to have in a graded algebra, while $x^4 + x = 0$ is not. In a graded algebra, the subspaces

$$\mathbb{A}_k = \text{span of monomials in generators of total degree } k \tag{45}$$

intersect only in zero for different $k$. Mathematicians put a circle around the plus sign in

$$\mathbb{A} = \bigoplus_k \mathbb{A}_k \tag{46}$$

to stress this fact. One says that (46) is a *direct sum*.

Each $\mathbb{A}_k$ is a finite-dimensional linear space over $\mathbb{F}$ and its dimension $\dim_\mathbb{F} \mathbb{A}_k$ is a number which may be an interesting combinatorial function of $k$. For a combinatorial classic, consider the example

$$\mathbb{A}_\mathbb{W} = \mathbb{Q}[x_1, x_2, x_3, x_4, x_5, x_6]/(x_1^{m_1+1} = 0, \ldots, x_6^{m_6+1} = 0),$$

where the degrees of variables, which we write as a vector, are the denominations of the Korean won coins

$$\deg x = (1, 5, 10, 50, 100, 500).$$

The reader should check that

$$\dim \mathbb{A}_{\mathbb{W},k} = \text{number of ways to pay } k \text{ won}$$

$$\text{using } \leq m_i \text{ coins of each denomination}.$$

For instance, if all $m_i = 1$, we get the sequence

| $k$ | 0 | 1 | 2 | 3 | 4 | 5 | 6 | 7 | 8 | 9 | 10 | 11 | 12 | ... |
|---|---|---|---|---|---|---|---|---|---|---|---|---|---|---|
| $\dim \mathbb{A}_{\mathbb{W},k}$ | 1 | 1 | 0 | 0 | 0 | 1 | 1 | 0 | 0 | 0 | 1 | 1 | 0 | ... |

, (47)

while for all sufficiently large $m_i$ one gets the sequence

| $k$ | 0 | 1 | 2 | 3 | 4 | 5 | 6 | 7 | 8 | 9 | 10 | 11 | 12 | ... |
|---|---|---|---|---|---|---|---|---|---|---|---|---|---|---|
| $\dim \mathbb{A}_{\mathbb{W},k}$ | 1 | 1 | 1 | 1 | 1 | 2 | 2 | 2 | 2 | 2 | 4 | 4 | 4 | ... |

, (48)

The number 4 here comes from the fact that $10 = 5 + 5 = 5 + 1 + \cdots + 1 = 1 + \cdots + 1$ are all valid ways to pay 10 won.

The reader should further check that when all $m_i$ are finite, the sequence $\dim \mathbb{A}_{\mathbb{W},k}$ is always palindromic, that is, equals to itself read backwards. Equivalently,

$$\dim \mathbb{A}_{\mathbb{W},k} = \dim \mathbb{A}_{\mathbb{W}, \text{topdeg}-k}, \quad \text{topdeg} = \sum_{i=1}^{6} m_i \deg x_i, \tag{49}$$

where the top degree is the total sum of money in our possession[8]. In particular, the sequence starts and ends with

$$\dim \mathbb{A}_{\mathbb{W},0} = \dim \mathbb{A}_{\mathbb{W}, \text{topdeg}} = 1,$$

---

[8] Instead of paying $k$ won, we can just give all our money and ask for topdeg $-k$ won in change.

but in the middle can have many ups and down, as exemplified by (47). In particular, it is not in general a *unimodal* sequence, where unimodality is the concept we recall from the discussion following (10).

### 6.3. Hard Lefschetz property

Remarkably, in a closely related situation, one finds a sequence which is not just palindromic, but also unimodal. Every country in the Euro-area has its own euro cent coin, so there are many different coins each worth €0.01. The number of ways to pay $k$ cents, using at most $m_i$ cents from the country number $i = 1, \ldots, N$, is related to the algebra

$$\mathbb{A}_\mathbb{E} = \mathbb{Q}[x_1, \ldots, x_N]/(x_1^{m_1+1} = 0, \ldots, x_N^{m_N+1} = 0),$$

with

$$\deg x = (1, \ldots, 1).$$

If $m_i = 1$ for all $i$ then $\dim \mathbb{A}_{\mathbb{E},k} = \binom{N}{k}$ is the binomial coefficient. If all $m_i \geq 2$ then

$$\dim \mathbb{A}_{\mathbb{E},k} = 1, N, \frac{N(N+1)}{2}, \ldots, \frac{N(N+1)}{2}, N, 1,$$

where the last 1 occurs in the maximal degree $\sum m_i$. Moreover, for any $m_i$'s, it will be always a unimodal sequence! (This is a combinatorial classic, for which the reader can try finding her or his own proof. It is probably the easiest to prove that the sequence of logarithms $\ln \dim \mathbb{A}_{\mathbb{E},k}$ is a *concave* function of $k$; see for example [11].)

There is certain algebraic property of $\mathbb{A}_\mathbb{E}$ that is stronger than the symmetry (49) and the unimodality. It concern multiplication by the element

$$\omega = \sum_{i=1}^{N} x_i \in \mathbb{A}_{\mathbb{E},1}.$$

The following diagram, plotted for $(m_1, m_2, \ldots) = (3, 1, 0, 0, \ldots)$ may help the reader visualize what this multiplication operator does.

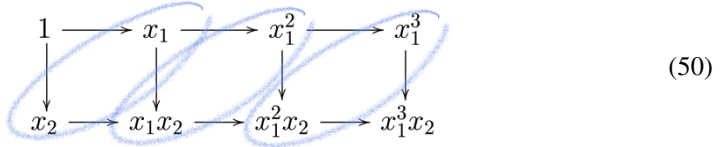

(50)

In (50), we have listed all monomials that do not reduce to zero on applying the rules $x_1^4 = x_2^2 = x_3 = \cdots = 0$. The graded pieces $\mathbb{A}_{\mathbb{E},k}$, $k = 0, 1, 2, 3, 4$ correspond to the diagonals in (50). The operations of multiplication by $x_1$ and $x_2$, when nonzero, are represented by the horizontal and vertical arrows, respectively. Thus, multiplication by $\omega$ is the sum of all outgoing arrows. We hope the reader has no problem visualizing the general case from this small example.

With these preparations, consider the multiplication map

$$\mathbb{A}_{\mathbb{E},i} \xrightarrow{\omega^{\text{topdeg}-2i}} \mathbb{A}_{\mathbb{E},\text{topdeg}-i}, \quad i < \tfrac{\text{topdeg}}{2}. \tag{51}$$

In the example (50), the possibilities for $i$ are $i = 0, 1$. For $i = 0$, we get

$$\omega^4 \cdot 1 = 4 x_1^3 x_2.$$

So, multiplication by $\omega^4$, identifies $\mathbb{A}_{\in,0} = \mathbb{Q}$ and $\mathbb{A}_{\in,4} = \mathbb{Q} \cdot x_1^3 x_2$. Similarly,

$$\omega^2 \cdot (c_1 x_1 + c_2 x_2) = c_1 x_1^3 + (c_2 + 2c_1) x_1^2 x_2, \quad c_1, c_2 \in \mathbb{Q},$$

which is easily seen to be injective and surjective (one is enough, since this is linear map between spaces of the same dimension). Thus, in the example (50), we observe that the maps (51) are isomorphisms[9].

In fact, the maps (51) are isomorphisms for the algebra $\mathbb{A}_\in$ for any $(m_1, m_2, \dots)$. The reader who wants to prove this directly will probably find it a good challenge.

In general, one says that a graded algebra $\mathbb{A}$ satisfies the hard Lefschetz property (HLP) if the multiplication maps (51) are isomorphisms for some $\omega \in \mathbb{A}_1$. Among other things, HLP implies that the multiplication map

$$\mathbb{A}_{\in,i} \xrightarrow{\omega} \mathbb{A}_{\in,i+1}, \quad i < \tfrac{\text{topdeg}}{2}, \tag{52}$$

is injective, whence the unimodality of the sequence $\dim \mathbb{A}_{\in,i}$.

Solomon Lefschetz was a very famous topologist, who proved[10] HLP for *cohomology* algebras of a certain class of *manifolds*; see Section 8. For the algebra $\mathbb{A}_\in$, the manifold in question is the product

$$\mathbb{P}^{m_1}(\mathbb{C}) \times \mathbb{P}^{m_2}(\mathbb{C}) \times \dots \mathbb{P}^{m_N}(\mathbb{C})$$

of complex projective spaces.

The geometric story of HLP and its generalizations is an immensely beautiful subject, about which we will try to say a few words in Section 8. There are two reasons our contact with this subject will be only tangential. First, going any deeper into this story requires a level of mathematical sophistication that is well beyond the style of these notes. Second, and more importantly, the work of [10] completely *bypasses* the old HLP story, creating a totally new combinatorial alternative for it. It is true that the old HLP served as in important inspiration and, in fact, the original proof of Theorem 1 relied on it. But progress in mathematics also sometime includes letting go of very beautiful constructions that are no longer logically required.

### 6.4. The graded Möbius algebra, finally

Given a matroid $M$, its graded Möbius algebra $\mathbb{H}(M)$ is defined as follows. As a linear space over $\mathbb{Q}$, it has a basis $y_F$ indexed by the flats $F$ of $M$. It is graded by

$$\deg y_F = \text{rank}(F).$$

---

[9] Mathematicians call a bijective map between two sets an isomorphism when it preserves some further structures that these sets possess. Multiplication by $\omega$ preserves multiplication by $\mathbb{A}_0$ and addition.

[10] Lefschetz's arguments were not entirely rigorous. Different correct proofs of HLP were given, in various geometric contexts, by Hodge, Chern, Deligne, and others. The influence of Lefschetz's work, however, was such that no one considers not naming this property after him.

Therefore,
$$|\mathcal{F}_r| = \dim \mathbb{H}(M)_r, \quad r = 0, 1, 2, \ldots.$$

The multiplication is defined by

$$y_F \, y_{F'} = \begin{cases} y_{F \vee F'}, & \text{rank}(F \vee F') = \text{rank}(F) + \text{rank}(F'), \\ 0, & \text{otherwise}, \end{cases} \tag{53}$$

where $F \vee F'$ is the minimal flat that contains $F$ and $F'$. Recall we have used this notation to denote the line $P_1 \vee P_2$ spanned by points $P_1$ and $P_2$, etc.

In particular, $y_\emptyset \in \mathbb{H}(M)_0$ is the identity for this product. Also note that

$$y_F \cdot \mathbb{H}(M) \subset \text{span}\left(\{y_{F'}\}_{\text{such that } F \subset F'}\right). \tag{54}$$

Theorem 2 easily follows from the following property of multiplication by the element[11]

$$\omega = \sum_{F \in \mathcal{F}_1} y_F \tag{55}$$

in the algebra $\mathbb{H}(M)$.

**Theorem 3** ([10]). *For $r \leq r'$ and $r + r' \leq \text{rank}(M)$ the linear map*

$$\mathbb{H}(M)_r \xrightarrow{\text{multiplication by } \omega^{r'-r}} \mathbb{H}(M)_{r'} \tag{56}$$

*is injective.*

Here is how Theorem 3 implies Teorem 2. Let

$$\mathsf{A} = (a_{F',F})$$

be the matrix of multiplication by $\omega^{r'-r}$ in the bases $\{y_F\} \subset \mathbb{H}(M)_r$ and $\{y_{F'}\} \subset \mathbb{H}(M)_{r'}$. See Appendix A for a reminder about matrices and bases.

By (54), the matrix entry $a_{F',F}$ vanishes unless $F \subset F'$. Since $\mathsf{A}$ is injective, it has is a square invertible submatrix $\mathsf{A}'$ of size $|\mathcal{F}_r|$. Since $\mathsf{A}'$ is invertible, its determinant

$$\det \mathsf{A}' \neq 0 \tag{57}$$

is not zero. Since $\det \mathsf{A}' \neq 0$, there is at least one nonzero term in the formula (110). The corresponding permutation determines an injective matching of flats in $\mathcal{F}_r$ to flats in $\mathcal{F}_{r'}$. Quod erat demonstrandum.

Note how the logic of the proof goes from combinatorics to linear algebra and back. A linear map takes a basis vector to a linear combination of basic vectors, and this gives linear maps important extra flexibility. The argument above is about how one can go back, and obtain an injective maps between bases from an injective linear map.

---

11    More generally, one can replace each $y_F$ in (55) by $c_F y_F$, where $c_F$ is an arbitrary positive rational number.

In broadest possible strokes, the strategy of the proof of Theorem 3 is the following. The authors of [10] construct a larger graded linear space

$$\mathbb{H}(M)_r \subset \mathbb{IH}(M)_r, \quad r = 0, 1, \ldots, \text{rank}(M),$$

which is no longer an algebra, but still has multiplication by elements of $\mathbb{H}(M)$. Mathematicians say $\mathbb{IH}(M)$ is a *module* for the algebra $\mathbb{H}(M)$. Since it is a module, it makes sense to consider the map

$$\mathbb{IH}(M)_r \xrightarrow{\text{multiplication by } \omega^{\text{rank}(M)-2r}} \mathbb{IH}(M)_{\text{rank}(M)-r}, \tag{58}$$

for $r < \frac{1}{2}\text{rank}(M)$. Evidently, (58) being an isomorphism implies that (64) is injective. It is this HLP for the map (58) that is really the key to Theorems 1, 2, and 3.

## 7. The big induction

In a computer code, it is sometimes very convenient to allow a procedure to call another instance of itself. Of course, if done carelessly, this can easily lead to an infinite loop and failure. To make sure the code terminates, there have to be, first, some base cases, when the procedure returns the answer without calling itself, and, second, it should each time call itself on smaller or simpler input, which gets closer and closer to a base case.

Imagine we already coded a data type `matroid` and we want to code, in some imaginary relative of the C programming language, a procedure

```
 print_theorem_3_proof(matroid M){
     if (rank M < 2){
         print("multiplication by ω⁰ is an isomorphism")
         ...           ,
```
(59)

where we have already indicated the base case. Definitely, a matroid $M'$ is simpler than $M$ if that has fewer points, so it is OK for this procedure to call `print_theorem_3_proof(M')` inside itself for such $M'$.

There are two important ways to construct a smaller matroid from $M$ and a flat $F$ of $M$. They are denoted $M^F$ and $M_F$. The matroid $M^F$ keeps only points $P_i$ and flats $F'$ contained in $F$. The matroid $M_F$ keeps only those flats $F'$ that contain $F$. The latter are determined by which points we should add to $F$ to get $F'$, hence $M_F$ is a matroid on the points $P_j$ that are *not* contained in $F$.

Of course, to have a mathematical proof of Theorem 3 it is not necessary to actually run the procedure. It is enough to know that a proof for $M$ can be found if we have a proof for all smaller matroids $M_F$ and $M^F$, and in the base case. Mathematicians call such proofs a proof by *induction*.

A very important insight from [10] is that it is much more natural to prove a *stronger* theorem than Theorem 3. In an inductive proof, there is always a tension between trying to prove too much or too little. The logic of induction says that we can get from the result for $M'$ to the result for $M$. So, assuming we can prove the statement for $M'$, we can prove it for $M$. There is a certain climb between $M'$ and $M$, and it becomes impossible if the starting point is too low or the goal is too high.

Analogies aside, what the authors of [10] actually prove is the whole *Kähler package* for the space $\mathbb{IH}(M)$. In addition to HLP, this package includes a nondegenerate bilinear form

$$(\,\cdot\,,\,\cdot\,) : \mathbb{IH}(M)_i \times \mathbb{IH}(M)_{\text{rank}-i} \to \mathbb{Q}, \tag{60}$$

such that for $\alpha \in \mathbb{IH}(M)_i$ we have

$$\omega^{\text{rank}-2i+1}\alpha = 0 \quad \Rightarrow \quad (-1)^i(\omega^{\text{rank}-2i}\alpha, \alpha) > 0. \tag{61}$$

For realizable matroids, these properties have an interpretation and history in topology, at which we will hint in Section 8. Namely (60) is the Poincaré duality and (61) are the Hodge-Riemann relations. But as we have already stressed at many points of this narrative, the amazing feature of Theorem 3 is that it works with no input from topology or algebraic geometry, and applies to absolutely all matroids, realizable or not.

The body of the procedure (59) is a marvel of combinatorics and combinatorial algebra, and it is way beyond the sophistication level of these notes to try to look any further in it. Let's just say it is not at all simple. There is a reason mathematics like this is recognized by the highest prizes in mathematics. In fact, it is miracle that some people can construct proofs like this.

An interested reader will find further reading suggestions in Section 9. We should also mention that Theorem 3 is not first time a combinatorial replacement of Hodge theory appears in mathematics.

June Huh says: "*Important precursors include the intersection cohomology $\mathbb{IH}(P)$ of a convex polytope $P$ [29] and the Soergel bimodule $\mathbb{IH}(w)$ for a Coxeter group element $w$ [17]. Both $\mathbb{IH}(P)$ and $\mathbb{IH}(w)$ satisfy Poincaré duality, the hard Lefschetz theorem, and the Hodge–Riemann relations, and these reveal fundamental properties of $P$ and $w$: The generalized lower bound conjecture for the number of faces in the case of $P$ [29, 41] and the nonnegativity conjecture for the coefficients of Kazhdan–Lusztig polynomials of Bruhat intervals in the case of $w$ [17, 31]. Each of the known proofs of the three combinatorial Kähler packages involves numerous details that are unique to that specific case.*"

Speaking of Kazhdan–Lusztig polynomials, the authors of [10] prove, in fact, much more than we managed to explain in these notes. In particular, they prove the nonnegativity of KL polynomials for all matroids.

I hope the readers share the narrator's sense of awe at this absolutely amazing mathematics and join me in warmest congratulations on it being recognized by the Fields Medal. I also hope the readers got the sense that today's mathematics is not just extraordinarily powerful, but also concrete, understandable, and fun, once one finds the right idea and the right point of view. While finding that right point of view is not at all easy, my biggest hope is to have inspired my youngest readers to believe that mathematics can be beautiful and rewarding, both as a subject and as a profession. Maybe this is also a good place for me to thank June Huh and Gil Kalai for this special opportunity to be introduced to their wonderful subject.

## 8. Inspirations from topology
### 8.1. Cohomology

Consider a graph $\Gamma$ drawn on a torus $\Sigma$, that is, on the surface of a bagel. By definition, a graph is a collection of *vertices* and *edges*. Since it is drawn on a surface $\Sigma$, it partitions $\Sigma$ into some regions that we will call *faces*. The vertices, edges, and faces are the 0-, 1-, and 2-dimensional objects in Figure (62).

We assume, and it is an important assumption, that every face is a polygon. For instance, the unique face in Figure (62) is obtained by gluing the opposite sides of the hexagon in (65). For a very different example of a graph on the torus, the reader may take the fine square mesh representing the torus in Figure (62)

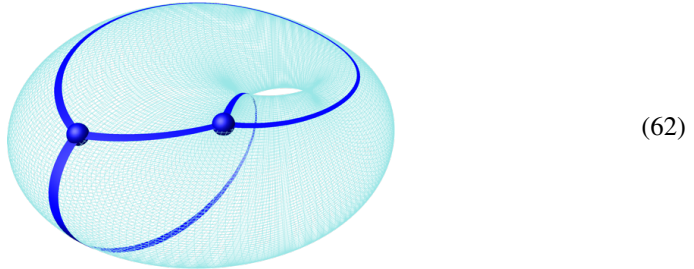

(62)

Let $f_0$ be a function defined on the vertices $\{V_i\}$ of $\Gamma$. We define its *gradient* $df_0$ as follows. This will be a function of an *oriented* edge $E_{ij}$ of $\Gamma$. If $E_{ij}$ goes from the vertex $V_i$ to the vertex $V_j$, schematically

$$V_i \bullet \xrightarrow{E_{ij}} \bullet V_j,$$

then

$$df_0(E_{ij}) = f(V_j) - f(V_i). \tag{63}$$

For the opposite orientation of the edge, one gets the opposite sign:

$$df_0(\overleftarrow{E}) = -df_0(\overrightarrow{E}). \tag{64}$$

This elementary construction is found in abundance is both theoretical and applied situations. For instance one may interpret $df_0$ as the current through edges of $\Gamma$ generated by a potential function $f_0$ defined on its vertices.

Now let $f_1$ be a function on oriented edges satisfying the sign rule (64). We may interpret $f_1$ as a flow or a vector field going along the edges of $\Gamma$. Given on oriented face $F$, its *boundary*

$$\partial F = E_{ij} \cup E_{jk} \cup \ldots$$

is a collection of edges, which gets an orientation from the orientation of $F$; see Figure 65.

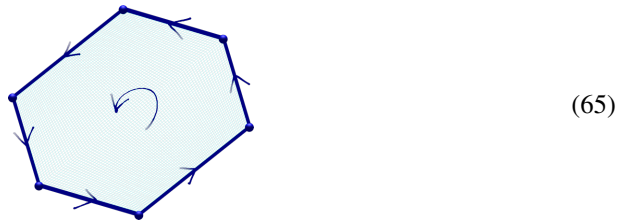

(65)

We define
$$df_1(F) = \sum_{E \in \partial F} f_1(E) . \tag{66}$$

This has a natural interpretation as the circulation, or the curl, of the flow $f_1$ around the face $F$. This again changes sign upon switching the orientation of $F$. The reader may want to pause and write (63) in a form that resembles (66).

So far, we did not specify the range of the functions $f_i$. Let them take values in some field $\mathbb{F}$. Functions on vertices form a linear space over $\mathbb{F}$ that we will denote $\Omega^0(\Gamma)$, and similarly for functions $\Omega^i(\Gamma)$, $i = 1, 2$ on oriented edges and faces. We always assume these functions change sign as in (64) upon switching the orientation. To fix a basis in these spaces, we can fix some orientation of each edge and face arbitrarily.

Above, we have constructed two linear maps
$$\Omega^0(\Gamma) \xrightarrow{d_0} \Omega^1(\Gamma) \xrightarrow{d_1} \Omega^2(\Gamma) , \tag{67}$$

where we marked the two maps $d$ with lower indices for notational convenience. The key property of (67) is that the composition
$$d^2 = d_1 d_0 = 0 \tag{68}$$

is zero. This is known in may different guises, e.g. the circulation of a gradient vector field is zero, and reflects the geometric fact that $\partial^2 = 0$, meaning that a boundary has no boundary.

A classical question appearing in many branches of mathematics is: does every vector field with zero curl comes from a potential? In other words is it true that the kernel $\operatorname{Ker} d_1$ equals the image $\operatorname{Im} d_0$? Or, using the language introduced in Section A.4, is the sequence (67) *exact* in the middle term?

More generally, one calls a sequence of maps composing to zero like (67) a *complex*, with the stress on the second syllable. From $d^2 = 0$, we see that $\operatorname{Im} d_{i-1} \subset \operatorname{Ker} d_i$ and one defines the *cohomology* groups of a complex by
$$H^i = \operatorname{Ker} d_i \,/ \operatorname{Im} d_{i-1} .$$

In (67) and in general we assume that the maps $d_i$ not indicated are the zero maps. The image of a zero map is the zero subspace and the kernel of a zero map is the whole space.

Clearly,
$$\operatorname{Ker} d_0 = \text{constant functions} = \mathbb{F} ,$$

hence $\dim H^0 = 1$. It fun to check that
$$\operatorname{Im} d_1 = \operatorname{Ker} \left( \Omega^2(\Gamma) \xrightarrow{\int_\Sigma} \mathbb{F} \right)$$

where the integration map is defined by
$$\int_\Sigma f_2 = \sum_{\text{all faces } F} f_2(F) ,$$

with the orientation of each face induced by some chosen orientation of $\Sigma$. Therefore, $\dim H^2 = 1$. The remaining dimension $\dim H^1$ we can infer from

$$\dim H^0 - \dim H^1 + \dim H^2 = \dim \Omega^0 - \dim \Omega^1 + \dim \Omega^2$$
$$= |\text{vertices}| - |\text{edges}| + |\text{faces}|$$
$$= \text{Euler characteristic of } \Sigma$$
$$= 0, \qquad (69)$$

where the first equality is a general property of all complexes which follows from (98) and (100), and where at the last step we used Euler's formula, one of the first topological results in mathematics.

Thus $\dim H^1 = 2$ and, in fact, a vector field $f_1$ is a gradient if and only if its circulation around each face *and* around the two blue loops as in Figure (70) vanishes.

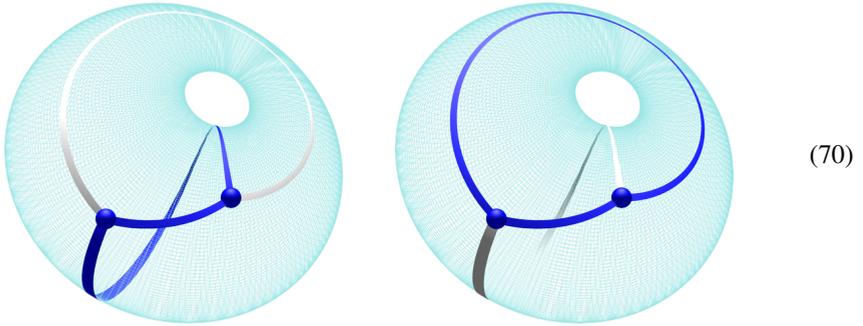

(70)

While we started with a graph on a torus, the eventual outcome of our computation is really about the torus itself and not about any particular graph drawn on it. One way to see it, is to take a a refinement $\Gamma'$ of the graph $\Gamma$. This means that every edge of $\Gamma$

$$E_{ij} = E'_{ik_1} \cup E'_{k_1 k_2} \cup \cdots \cup E'_{k_l j}$$

is a union of edges of $\Gamma'$. It follows that every face of $\Gamma$ is a union of faces of $\Gamma'$. A refinement of a face may look something like Figure (71)

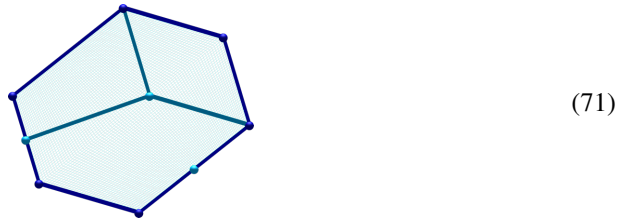

(71)

From the flows and circulations in the graph $\Gamma'$ we can compute the flows and circulations in the graph $\Gamma$. This gives vertical maps in the diagram

$$\begin{array}{ccccc} \Omega^0(\Gamma) & \xrightarrow{d_0} & \Omega^1(\Gamma) & \xrightarrow{d_1} & \Omega^2(\Gamma) \\ \uparrow & & \uparrow & & \uparrow \\ \Omega^0(\Gamma') & \xrightarrow{d'_0} & \Omega^1(\Gamma') & \xrightarrow{d'_1} & \Omega^2(\Gamma') \end{array} \qquad (72)$$

which can be seen to identify the cohomology. Since any two graphs $\Gamma_1$ and $\Gamma_2$ have a common refinement $\Gamma'$, the cohomology is really attached to the torus $\Sigma$ and not to any particular decomposition of $\Sigma$ into vertices, edges, and faces. What is really fundamental about the torus is that all possible loops in $\Sigma$, up to boundaries, span a two-dimensional space with a basis plotted in Figure (70).

We hope it is easy for the reader to imagine the generalization of this story in which one replaces the torus $\Sigma$ by any *topological space* that can be glued out of polytopes of different dimension. The elementary story told above was known in the XIX century, and since then topologists have developed really powerful tools to attach various topological invariants, including cohomology $H^i(X, \mathbb{F})$, to topological spaces $X$.

The torus is not only a topological space, it is also a *complex algebraic variety*. Namely, consider the solutions

$$\{P(x_0, x_1, x_2) = 0\} \subset \mathbb{P}^2(\mathbb{C}), \tag{73}$$

where $P$ is a homogeneous polynomial of degree 3, meaning that

$$P(tx_0, tx_1, tx_2) = t^3 P(x_0, x_1, x_2), \quad \text{for any } t \in \mathbb{C}. \tag{74}$$

While the coordinates $(x_0, x_1, x_2)$ on $\mathbb{P}^2$ are defined only up to proportionality, by (74) the zero set (73) is defined unambiguously. If the partial derivatives $\frac{\partial}{\partial x_i} P$ do not vanish simultaneously, then (73) is a torus. Not to be outdone by the topologists, algebraic geometers have defined equally powerful cohomology theories for algebraic varieties. These agree with the topological definitions over the field $\mathbb{C}$ of complex numbers.

It is a really inspiring lesson in the unity of mathematics that different cohomology theories, with very different starting points and emphasis on very different geometric objects, in the end all agree on their common domains of applicability.

### 8.2. Multiplication and Poincare duality

We were interested in cohomology because of the graded algebra structure on the direct sum

$$H^\bullet(X, \mathbb{F}) = \bigoplus H^i(X, \mathbb{F}),$$

that is, because of the multiplication operation

$$H^i(X, \mathbb{F}) \times H^j(X, \mathbb{F}) \xrightarrow{\cup} H^{i+j}(X, \mathbb{F}). \tag{75}$$

The product (75) exists for *very* abstract reasons. For any topological space $X$ there is the diagonal map

$$X \to X \times X,$$

sending a point $x$ to the pair $(x, x)$. A map between topological spaces induces a map on cohomology the other way. Using the Künneth isomorphism

$$H^\bullet(X \times X, \mathbb{F}) = H^\bullet(X, \mathbb{F}) \otimes H^\bullet(X, \mathbb{F}), \tag{76}$$

where $\otimes$ denotes the tensor product, one obtains (75). A less general but more intuitive description says that (75) is dual to *intersection*, and it goes as follows.

By definition, the complex of *dual* maps

$$\Omega^0(\Gamma)^\vee \xleftarrow{d_0^\vee} \Omega^1(\Gamma)^\vee \xleftarrow{d_1^\vee} \Omega^2(\Gamma)^\vee\,, \qquad (77)$$

computes the *homology* groups $H_i(\Sigma, \mathbb{F})$. The geometric meaning of (77) is transparent. The bases in $\Omega^i(\Gamma)^\vee$, $i = 0, 1, 2$, are indexed by vertices, edges, and faces, respectively, and the maps are the boundary maps $\partial$. Collectively, the vertices, edges, and faces are called *cells*, and $\mathbb{F}$-linear combinations of cells are called *chains*. Chains with no boundary are called *cycles*. Thus, homology describes cycles up to boundaries.

For the torus $\Sigma$ the homology complex (77) can be identified with the cohomology complex for the *dual graph* $\Gamma^\vee$, where the dual graph to (62) can be seen in (78).

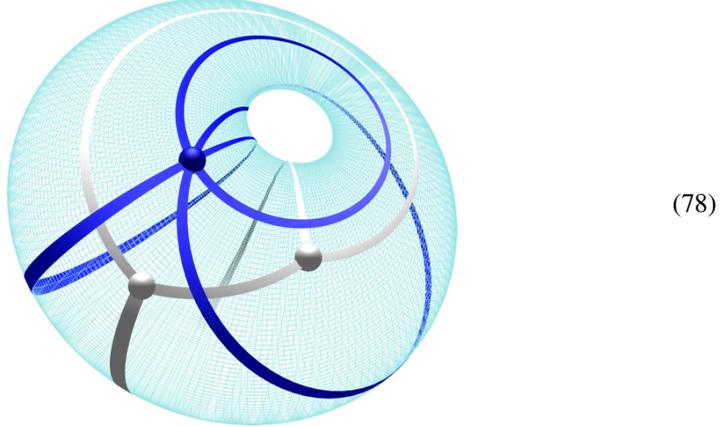

(78)

The vertices, edges, and faces of the dual graph correspond to the faces, edges, and the vertices of the original graph, respectively. Moreover, each cell intersects the dual cells in exactly one point. This gives the Poincaré duality isomorphism

$$H_i(\Sigma, \mathbb{F}) \cong H^{\dim \Sigma - i}(\Sigma, \mathbb{F})\,. \qquad (79)$$

It works just the same for a doughnut with $g = 2, 3, \dots$ holes and for any oriented closed (meaning, compact and without boundary) manifold $M$.

Manifolds are particularly nice topological spaces that, in a certain technical sense, look just like the linear space in the vicinity of every point. The linear space $\mathbb{R}^n$ is a manifold, but not a compact manifold. The $n$-dimensional sphere $S^n$ and also the real and complex projective spaces are closed manifolds. Recall we insisted that $\frac{\partial}{\partial x_i} P \ne 0$ for some $i$ at every point of (73). This was to make sure that (73) defines a manifold.

For cycles, one would like to define an intersection product

$$H_i(X, \mathbb{F}) \times H_j(X, \mathbb{F}) \xrightarrow{\cap} H_{i+j-\dim X}(X, \mathbb{F})\,, \qquad (80)$$

that would turn into the $\cup$-product upon the identification (79). It doesn't really make sense for a general topological space, since it is not even clear what notion of dimension one should use in (80). But on a manifold, it works beautifully, especially if one intersects cycles defined using a graph $\Gamma$ with the cycles defined using the *dual* graph $\Gamma^\vee$. (Recall that any two graphs $\Gamma$ and $\Gamma'$ define the *same* space of cycles up to boundaries.)

One very important detail here that one should count intersections with a sign according to orientation. This is crucial to make boundaries have zero intersection with any cycle. Lets look at the intersection of a cycle $\gamma$ with a boundary $\partial F$ in Figure (81). If we keep track of the orientations, we can tell whether $\gamma$ enters or exists $F$ at a given point of intersection. Hence if we count intersections with signs then we get $\gamma \cap \partial F = 0$.

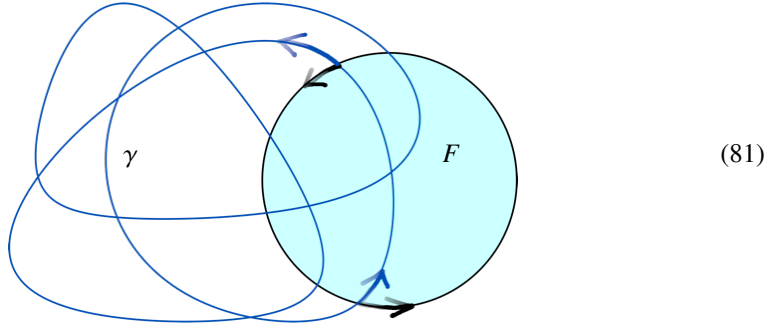 (81)

As a result, cohomology $H^\bullet(\Sigma, \mathbb{F})$ not commutative but rather *supercommutative*[12]

$$\alpha_1 \cup \alpha_2 = (-1)^{d_1 d_2} \alpha_2 \cup \alpha_1, \quad \alpha_i \in H^{d_i}. \tag{82}$$

If the reader has not seen it before, it is a good exercise to work out multiplication for $H^\bullet(\Sigma, \mathbb{F})$ and also for $H^\bullet(\Sigma_g, \mathbb{F})$, where $\Sigma_g$ is the surface of the doughnut with $g$ holes.

While the signs in (82) are an important fact of nature, they are of little concern for us here since we are interested in even cohomology, which is commutative. In particular, topological intersection of algebraic cycles, that is, those defined by polynomial equations, is commutative and agrees with its counterpart in algebraic geometry.

For cycles of complementary dimension, we can interpret interpret the isomorphism (79) as the Poincaré duality pairing:

$$(\cdot, \cdot) : H^i(X, \mathbb{F}) \times H^{\dim X - i}(X, \mathbb{F}) \to \mathbb{F}. \tag{83}$$

The purely combinatorial relative of this pairing appears in (60). One can also interpret (83) as the composition of the cup product with the isomorphism $H^{\dim X}(X, \mathbb{F}) \cong H_0(X, \mathbb{F}) = \mathbb{F}$.

It is an excellent student project to prove that

$$H^\bullet(\mathbb{P}^d(\mathbb{C}), \mathbb{F}) = \mathbb{F}[x]/(x^{d+1} = 0), \quad x \in H^2, \tag{84}$$

where $x$ is dual to the class of the hyperplane in $H_{2d-2}(\mathbb{P}^d(\mathbb{C}))$. This formalizes the table (44). Similarly,

$$\mathbb{A}_{\in, k} = H^{2k}\left(\mathbb{P}^{m_1}(\mathbb{C}) \times \cdots \times \mathbb{P}^{m_N}(\mathbb{C}), \mathbb{Q}\right). \tag{85}$$

### 8.3. The Hard Lefschetz property

The title of this subsection prompts the question: for which even-dimensional manifolds $X$ is there a class $\omega \in H^2(X, \mathbb{Q})$ such that the multiplication map

$$H^i(X, \mathbb{Q}) \xrightarrow{\omega^{\frac{1}{2} \dim_\mathbb{R} X - i}} H^{\dim_\mathbb{R} X - i}(X, \mathbb{Q}), \quad i < \tfrac{1}{2} \dim_\mathbb{R} X, \tag{86}$$

---

**12** The tensor product in Künneth theorem (76) should also be understood with signs.

is an isomorphism? Here by the dimension of $X$ we mean its real dimension, even though $X$ may have been originally defined as a complex manifold. Thus, the $\dim_{\mathbb{R}} \Sigma = 2$ for the torus (73). It is not enough for $X$ to be smooth: the even dimensional spheres $X = S^{2k}$ have $H^2(S^{2k}) = 0$, and hence no chance of satisfying (86) for $k > 1$.

One classical answer is that a smooth projective $X \subset \mathbb{P}^N(\mathbb{C})$ satisfied the Hard Lefschetz property (86) with $\omega$ dual to the class of a hyperplane section. Projective means that $X$ is defined by polynomial equations just like (73), and smooth means it is a manifold. A more general class of *Kähler* manifold also satisfies the HLP.

If $X \subset \mathbb{P}^N(\mathbb{C})$ is not smooth then it is called *singular*. For singular $X$, there is a more delicate cohomology theory that satisfies the HLP. It is called intersection cohomology and its developments is one of the true highlights of the geometry and topology, achieved in the 1970's and 1980's by Mark Goresky, Robert MacPherson, Pierre Deligne, Alexander Beilinson, Joseph Bernstein, and other amazing mathematicians.

For a matroid linearly realisable over $\mathbb{C}$, the Möbius algebra $\mathbb{H}(M)$ is the cohomology algebra of the variety $\overline{Y}$ associated to $M$ in (38). Note that the generators $\mathbb{H}(M)$ commute and square to zero, so it is a quotient of $H^\bullet\left((\mathbb{P}^1)^n\right)$. The module $\mathbb{IH}(M)$ is the intersection cohomology of the same variety $\overline{Y}$. Thus, there is a topological proof the HLP for $\mathbb{IH}(M)$ for realizable matroids. This was used in the original proof of Theorem 1 given in [22].

Hard Lefschetz property for cohomology and intersection cohomology has a history of very powerful application to combinatorial problems. One great example is Richard Stanley's proof of McMullen's conjectural characterization of $f$-vectors[13] of simplicial convex polytopes. (Stanley proved the necessity of McMullen's conditions, the sufficiency was proven about the same time by Billera and Lee.) See [42] for a discussion of this and other combinatorial applications of the HLP.

## 9. Further reading

The *Quanta Magazine* has published popular accounts of these and related developments, see [13, 20].

Among surveys written by top experts in the field, one should mention [5, 28], including expositions by June Huh himself [23, 24].

Among textbook introductions to different areas mathematics discussed in our narrative, the reader will surely find something which suits her or his interests and style among [4, 21, 38, 40, 43].

I hope the reader has a lot of fun studying these sources as well as the original articles [10, 22].

---

**13**  For a convex polytope, its $f$-vector records the number of faces of each dimension.

## A. A rice bowl of linear algebra
### A.1. Linear equations

A system of $N$ linear equations for $M$ unknowns $x_1, \ldots, x_M$,

$$\begin{cases} a_{11}x_1 & + & \ldots & + & a_{1M}x_M & = c_1 \\ \vdots & & & & \vdots & \vdots \\ a_{N1}x_1 & + & \ldots & + & a_{NM}x_M & = c_N \end{cases} \tag{87}$$

is conveniently written in matrix notation

$$\mathsf{A}\mathsf{x} = \mathsf{c}, \tag{88}$$

where

$$\mathsf{A} = \begin{bmatrix} a_{11} & \ldots & a_{1M} \\ \vdots & & \vdots \\ a_{N1} & \ldots & a_{NM} \end{bmatrix}, \quad \mathsf{x} = \begin{bmatrix} x_1 \\ \vdots \\ x_M \end{bmatrix}, \quad \mathsf{c} = \begin{bmatrix} c_1 \\ \vdots \\ c_N \end{bmatrix}. \tag{89}$$

Solutions of (87) are unchanged, if we multiply $i$th equation, where $i = 1, \ldots, N$ by a nonzero number $t$, or add to the $i$th equation $t$ times the $j$th equation. Such transformation are called *elementary*. In matrix form, they have the form

$$(\mathsf{A}, \mathsf{c}) \to (\mathsf{g}_{ij}(t)\mathsf{A}, \mathsf{g}_{ij}(t)\mathsf{c}),$$

where $\mathsf{g}_{ij}(t)$ is an *elementary* matrix, which means a matrix of the following form

$$\underbrace{\begin{bmatrix} 1 & 0 & 0 & \ldots & 0 \\ 0 & 1 & 0 & \ldots & 0 \\ 0 & 0 & 1 & & \vdots \\ \vdots & \vdots & & \ddots & 0 \\ 0 & 0 & \ldots & 0 & 1 \end{bmatrix}}_{\text{identity matrix}} \xrightarrow{\text{put } t \neq 0 \text{ in } i\text{th row and } j\text{th column}} \mathsf{g}_{ij}(t). \tag{90}$$

The rice in our rice bowl is the following statement, called row reduction, or Gaussian elimination. Any system of linear equation can be transformed by elementary operations[14] to a unique matrix of the schematic form

$$\mathsf{A}_{\text{rowred}} = \begin{bmatrix} 0 & 0 & 1 & * & 0 & 0 & * & 0 & * \\ 0 & 0 & 0 & 0 & 1 & 0 & * & 0 & * \\ 0 & 0 & 0 & 0 & 0 & 1 & * & 0 & * \\ 0 & 0 & 0 & 0 & 0 & 0 & 0 & 1 & * \\ 0 & 0 & 0 & 0 & 0 & 0 & 0 & 0 & 0 \end{bmatrix}, \tag{91}$$

---

[14] In practical implementations of row reduction, it is very convenient to permute equations. Abstractly, however, a permutation of two equations may be achieved by elementary tranformations, as the reader will easily check.

where stars stand for some unspecified numbers. The 1's in (91) have to be in different rows and columns. A star $*$ can follow a 1 in a row, but not if there is a 1 in the same column. The number of 1's in (91) is called the rank of A.

We have

$$\text{solutions}\left(A\,x = c\right) = \text{solutions}\left(A_{\text{rowred}}\,x = c'\right), \tag{92}$$

where $c'$ is the result of applying the sequence of elementary tranformations $g_{ij}(t)$ to the vector c. For a reduced matrix, the solutions can be described very easily.

The zero rows in (91) lead the equations of the form $0 = c'_i$. These have either no solutions if $c'_i \neq 0$ or can be discarded if $c'_i = 0$. Thus (88) has solutions if and only if c satisfies $N - \text{rank}(A)$ linear equations given by $c'_i = 0$, $i = \text{rank}(A) + 1, \ldots, N$.

After we have dealt with the zero rows in (91), the remaining equations may be solved uniquely for the variables $x_j$ that have a 1 in their columns. All other variables are free parameters. Thus, when exists, the solutions are parametrized by $M - \text{rank}(A)$ many free parameters.

Row reduction is fundamental. Everything else in this section is a topping.

### A.2. Linear maps

In Section A.1, we never specified the algebraic nature of the variables $x_j$ or the coefficients $a_{ij}$ and $c_i$. The reader may have assumed they are real or rational numbers. In fact, they can be taken to be elements in any field $\mathbb{F}$ without changing anything at all in the analysis.

Column vectors x of size $M$ with entries in $\mathbb{F}$ form the coordinate linear space $\mathbb{F}^M$ of dimension $M$. It has operations of addition and multiplication by elements $t \in \mathbb{F}$, both defined coordinate by coordinate. A map

$$A : \mathbb{F}^M \to \mathbb{F}^N \tag{93}$$

is said to be linear if it preserves these operations, that is,

$$A(x + x') = A(x) + A(x'), \quad A(tx) = tA(x).$$

The reader should check that such maps are precisely the ones given by a matrix multiplication, and hence we can write A x in place of A(x). From what we just learned about linear equations, it follows that:

- A is injective if and only if $\text{rank}(A) = M$.

- A is surjective if and only if $\text{rank}(A) = N$.

- A is bijective, or an isomorphism, or *invertible* if and only if $\text{rank}(A) = N = M$. In particular, there is no isomorphism $\mathbb{F}^M \to \mathbb{F}^N$ if $M \neq N$.

- any isomorphism $g : \mathbb{F}^M \to \mathbb{F}^M$ is a product of elementary matrices $g_{ij}(t)$.

It is very important to remember that linear spaces have a lot of nontrivial isomorphisms $g : \mathbb{F}^M \to \mathbb{F}^M$. These can be composed and inverted, thus form a group denoted $GL(M, \mathbb{F})$.

When we act by $g \in GL(M, \mathbb{F})$, we say that the we change the basis, or do a linear change of coordinates. While coordinates provide a very concrete and convenient description of geometric objects, a truly geometric construction should work equally well in any linear coordinates.

Formula (91) describes a standard form to which a matrix can be brought by post-composing with an isomorphism, that it, by a change of basis in the target of the map. We invite the reader to check that by an independent[15] change of basis in the source and the target, a matrix can be brought to a particularly simple form

$$A_{\text{rowcolred}} = \begin{bmatrix} 1 & 0 & 0 & 0 & 0 & 0 & 0 & 0 & 0 \\ 0 & 1 & 0 & 0 & 0 & 0 & 0 & 0 & 0 \\ 0 & 0 & 1 & 0 & 0 & 0 & 0 & 0 & 0 \\ 0 & 0 & 0 & 1 & 0 & 0 & 0 & 0 & 0 \\ 0 & 0 & 0 & 0 & 0 & 0 & 0 & 0 & 0 \end{bmatrix}. \tag{94}$$

The number of ones in (94) is still the rank of $A$.

This means that, in some coordinates, a linear map from one linear space to another just forgets some coordinates, and pads the remaining ones by zeros. Thus a linear map from one linear space to another may be pictured as follows:

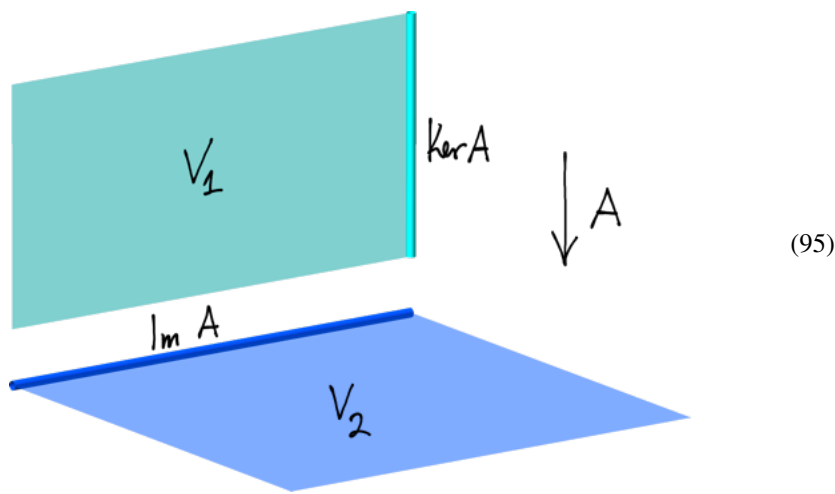

$$\tag{95}$$

The symbols $\operatorname{Ker} A$ and $\operatorname{Im} A$ will be explained in Section A.4.

### A.3. Abstract linear spaces

A set $V$ is called a linear space over a field $\mathbb{F}$ if it has a special element $0 \in V$, an operation of addition, and an operations of multiplication by $t \in \mathbb{F}$, satisfying the same rules as the corresponding operations in $\mathbb{F}^M$.

Any collection of vectors $v_1, \ldots, v_M$ determines a linear map

$$\mathbb{F}^M \to V$$

---

[15] Something much more interesting happens if the source and target are the *same* space and we have to use the *same* change of variables in both.

given by

$$\begin{bmatrix} x_1 \\ \vdots \\ x_M \end{bmatrix} \mapsto \sum x_i \mathsf{v}_i \,. \tag{96}$$

The linear space is said to have *finite dimension* if for some collection $\mathsf{v}_1, \ldots, \mathsf{v}_M$ the map (96) is surjective. By a version of row reduction, there is then a subset of $\mathsf{v}_i$'s for which the map (96) becomes an isomorphism. Such set of vectors is called a *basis* for $\mathsf{V}$, and its cardinality is denoted $\dim \mathsf{V}$. From what we already know, $\mathsf{V}$ has many bases, which can be all taken to each other by the group $GL(\mathsf{V}) = GL(\dim \mathsf{V}, \mathbb{F})$.

### A.4. Kernel, image, and quotient

Given a linear map

$$\mathsf{A} : \mathsf{V}_1 \to \mathsf{V}_2 \,,$$

one defines its *kernel* by

$$\operatorname{Ker} \mathsf{A} = \{\mathsf{v} \in \mathsf{V}_1, \mathsf{A}\mathsf{v} = 0\} \,. \tag{97}$$

This is a *linear subspace* of the source space of $\mathsf{A}$ of dimension

$$\dim \operatorname{Ker} \mathsf{A} = \dim \mathsf{V}_1 - \operatorname{rank} \mathsf{A} \,. \tag{98}$$

The projection in (95) is the projection *along* the kernel of $\mathsf{A}$. One defines the image of $\mathsf{A}$ as the image of this projection, or more formally

$$\operatorname{Im} \mathsf{A} = \{\mathsf{v}' \in \mathsf{V}_2, \mathsf{v}' = \mathsf{A}\mathsf{v} \text{ for some } \mathsf{v} \in \mathsf{V}_1\} \,. \tag{99}$$

This is a linear subspace in the target space of $\mathsf{A}$ of dimension

$$\dim \operatorname{Im} \mathsf{A} = \operatorname{rank} \mathsf{A} \,. \tag{100}$$

Thus, $\mathsf{A}$ may be factored as a projection and embedding

$$\mathsf{V}_1 \xrightarrow{\text{projection}} \operatorname{Im} \mathsf{A} \xrightarrow{\text{embedding}} \mathsf{V}_2 \,. \tag{101}$$

Mathematicians call a sequence of maps of the form

$$0 \to \operatorname{Ker} \mathsf{A} \xrightarrow{\text{embedding}} \mathsf{V}_1 \xrightarrow{\text{projection}} \operatorname{Im} \mathsf{A} \to 0 \tag{102}$$

a *short exact sequence*. This is an important word to remember and use. It is *exact* because the the kernel of each arrow in (102) is the image of the preceding map. It is short because it contains only 3 nontrivial terms. One also says that a projection is the *quotient* of $\mathsf{V}_1$ by its kernel.

### A.5. Dual vector spaces

For an $\mathbb{F}$-linear space $\mathsf{V}$, the space

$$\mathsf{V}^{\vee} = \{\text{linear functions } \mathsf{V} \xrightarrow{\xi} \mathbb{F}\}$$

is also an $\mathbb{F}$-linear space, since we can add linear functions and multiply them by numbers. The dual space to $\mathbb{F}^M$ is best visualized as the space of row vectors

$$(\mathbb{F}^M)^{\vee} = \{[\xi_1, \ldots, \xi_M]\}, \quad \xi_1, \ldots, \xi_M \in \mathbb{F}, \tag{103}$$

with

$$\xi(\mathsf{x}) = [\xi_1, \ldots, \xi_M] \begin{bmatrix} x_1 \\ \vdots \\ x_M \end{bmatrix} = \sum \xi_i x_i \stackrel{\text{def}}{=} \langle \xi, \mathsf{x} \rangle. \tag{104}$$

Here the notation $\langle \xi, \mathsf{x} \rangle$ is introduced to stress the symmetry between $\xi$ and $\mathsf{x}$. Mathematicians like to stress that, while $\mathsf{V}$ and $\mathsf{V}^{\vee}$ are vector spaces of the same dimension, there is no *natural*, that is, coordinate-independent identification between them. By constrast, the symmetry in (104) shows that $(\mathsf{V}^{\vee})^{\vee} = \mathsf{V}$ in a coordinate-free way for any finite-dimensional vector space.

For any linear subspace $\mathsf{V}' \subset \mathsf{V}$, there is the *annihilator* subspace

$$(\mathsf{V}')^{\perp} = \{\xi \text{ such that } \langle \xi, \mathsf{v}' \rangle = 0 \text{ for all } \mathsf{v}' \in \mathsf{V}'\}.$$

This is a subspace of $\mathsf{V}^{\vee}$ of dimension

$$\dim (\mathsf{V}')^{\perp} = \dim \mathsf{V} - \dim \mathsf{V}', \tag{105}$$

satisfying $((\mathsf{V}')^{\perp})^{\perp} = \mathsf{V}'$.

A linear map $\mathsf{A} : \mathsf{V}_1 \to \mathsf{V}_2$ induces the dual map

$$\mathsf{A}^{\vee} : \mathsf{V}_2^{\vee} \to \mathsf{V}_1^{\vee}$$

by precomposing a function with $\mathsf{A}$. In other words

$$\langle \mathsf{A}^{\vee} \xi, \mathsf{x} \rangle = \langle \xi, \mathsf{A} \mathsf{x} \rangle.$$

For row vectors as in (103), this is just *left* multiplication $\xi \mapsto \xi \mathsf{A}$. It is important that duality *reverses* the order of the composition

$$(\mathsf{A}_1 \mathsf{A}_2)^{\vee} = \mathsf{A}_2^{\vee} \mathsf{A}_1^{\vee}.$$

We have

$$\operatorname{Ker} \mathsf{A}^{\vee} = (\operatorname{Im} \mathsf{A})^{\perp}, \quad \operatorname{Im} \mathsf{A}^{\vee} = (\operatorname{Ker} \mathsf{A})^{\perp},$$

and in particular

$$\operatorname{rank} \mathsf{A}^{\vee} = \operatorname{rank} \mathsf{A}.$$

If we insist on identifying the row vectors with column vectors by switching the rows and columns then $\mathsf{A}^{\vee}$ becomes the *transposed matrix*

$$(a_{ij})^{\mathsf{T}} = (a_{ji}). \tag{106}$$

### A.6. Rank and rank

Let $P_1, \ldots, P_N$ be a collection of points in $\mathbb{P}^d(\mathbb{F})$, where the projective space is defined in Section 5.2. If

$$\mathbf{P} = N \times (d+1) \text{ matrix with rows } P_i,$$

then the equation

$$\mathbf{P}\mathbf{a} = 0, \quad \mathbf{a} = \begin{bmatrix} a_0 \\ \vdots \\ a_d \end{bmatrix}, \tag{107}$$

describes the hyperplanes containing the points $P_1, \ldots, P_N$. Thus

$$\begin{aligned} \mathrm{rank}(\{P_1, \ldots, P_N\}) &= \dim \mathrm{span}(\{P_1, \ldots, P_N\}) + 1 \\ &= (d+1) - \dim \text{ solutions of } (107) \\ &= \mathrm{rank}\,\mathbf{P}. \end{aligned} \tag{108}$$

More generally, the rank of any subset of $\{P_1, \ldots, P_N\}$ is the rank of the corresponding submatrix in $\mathbf{P}$.

For a practical computation of the rank, it is enough to bring a matrix by row operations to the row echelon form

$$A_{\text{row echelon}} = \begin{bmatrix} 0 & 0 & \star & * & * & * & * & * & * \\ 0 & 0 & 0 & 0 & \star & * & * & * & * \\ 0 & 0 & 0 & 0 & 0 & \star & * & * & * \\ 0 & 0 & 0 & 0 & 0 & 0 & 0 & \star & * \\ 0 & 0 & 0 & 0 & 0 & 0 & 0 & 0 & 0 \end{bmatrix}, \tag{109}$$

where $\star$ stands for some nonzero element of $\mathbb{F}$.

A theoretical formula for the rank may given using determinants; see Appendix B. Namely, the rank is the maximal size of a square submatrix with nonzero determinant.

### B. Determinant
### B.1. Formula

A matrix $A : \mathbb{F}^N \to \mathbb{F}^N$ is invertible if and only if $\det A \ne 0$, where $\det A$ is a certain polynomial in matrix elements $a_{ij}$. An explicit formula for this polynomial was needed in Section 6.4 to deduce the existence of a matching from equation (57).

This formula is a sum over *permutations* $\sigma$ of $\{1, \ldots, N\}$. It reads

$$\det A = \sum_{\text{permutations } \sigma} \mathrm{sgn}(\sigma)\, a_{1,\sigma(1)} a_{2,\sigma(2)} \cdots a_{N,\sigma(N)}, \tag{110}$$

where permutations and their signs are defined as follows. See further below for one possible explanation of the formula (110).

### B.2. Permutations

By definition, a permutation $\sigma$ of a finite set is a bijective map from a set to itself, like the following example

$$\begin{array}{ccccc} 1 & 2 & 3 & 4 & 5 \\ \downarrow & \downarrow & \downarrow & \downarrow & \downarrow \\ 1 & 2 & 3 & 4 & 5 \end{array} \tag{111}$$

for $N = 5$. A permutation has a sign defined by

$$\text{sgn}(\sigma) = (-1)^{|\text{crossings in (111)}|} = (-1)^{|\text{inversions}|}, \tag{112}$$

where an inversion of $\sigma$ is a pair $i < j$ such that $\sigma(i) > \sigma(j)$. For example, $(1, 2)$, $(1, 3)$, and $(4, 5)$ are the inversions for $\sigma$ in (111) and hence $\text{sgn}(\sigma) = -1$ for this particular permutation.

It is a nice exercise in the spirit of Figure (81) to check that

$$\text{sgn}(\sigma_1 \sigma_2) = \text{sgn}(\sigma_1) \, \text{sgn}(\sigma_2),$$

where $\sigma_1 \sigma_2$ denotes the composition of two permutations. In particular, the sign does not depend of how we order an $N$-element set. Compare the sign in (111) and below:

$$\begin{array}{ccccc} 1 & 2 & 4 & 3 & 5 \\ \downarrow & \downarrow & \downarrow & \downarrow & \downarrow \\ 1 & 2 & 4 & 3 & 5 \end{array}$$

### B.3. $N = 2$ case and the cohomology of the torus

Let

$$\mathsf{A} = \begin{bmatrix} a & b \\ c & d \end{bmatrix}$$

be a $2 \times 2$ matrix. Then

$$\det \mathsf{A} = ad - bc \tag{113}$$

and

$$\mathsf{A}^{-1} = \frac{1}{\det \mathsf{A}} \begin{bmatrix} d & -b \\ -c & a \end{bmatrix}, \tag{114}$$

over any field $\mathbb{F}$, as can be checked directly. Thus indeed we see that $\det \mathsf{A} \neq 0$ is equivalent to invertibility of $\mathsf{A}$.

Let us see what the cohomology of the torus $H^\bullet(\Sigma, \mathbb{F})$ can tell us about the formula (114). We hope the reader did the exercise suggested in Section 8.2 and remembers that

$$H^0(\Sigma, \mathbb{F}) = \mathbb{F}, \quad H^1(\Sigma, \mathbb{F}) = \mathbb{F}\gamma_1 \oplus \mathbb{F}\gamma_2, \quad H^2(\Sigma, \mathbb{F}) = \mathbb{F}\gamma_1 \cup \gamma_2, \tag{115}$$

for some basis $\{\gamma_1, \gamma_2\}$ of the 2-dimensional space $H^1(\Sigma, \mathbb{F})$. The description (115) means that $H^\bullet(\Sigma, \mathbb{F})$ is generated by $\gamma_1$ and $\gamma_2$ using the cup product, and the relations that these generators satisfy are

$$\gamma_1 \cup \gamma_1 = \gamma_2 \cup \gamma_2 = \gamma_1 \cup \gamma_2 + \gamma_2 \cup \gamma_1 = 0. \tag{116}$$

These can be written more compactly as follows

$$\text{for all } \gamma \in H^1, \quad \gamma \cup \gamma = 0. \tag{117}$$

Now, $H^1$ is 2-dimensional linear space over $\mathbb{F}$ with a basis, hence we can act by the matrix $\mathsf{A}$ in it. This action preserves the relation (117), and so induces the action on $H^2$. Since

$$(a\gamma_1 + c\gamma_2) \cup (b\gamma_1 + d\gamma_2) = (ab - cd)\,\gamma_1 \cup \gamma_2,$$

we conclude that $\mathsf{A}$ acts on $H^2$ as follows

$$H^2(\Sigma, \mathbb{F}) \xrightarrow{\text{multiplication by } \det \mathsf{A}} H^2(\Sigma, \mathbb{F}).$$

Writing the cup product as the Poincaré duality pairing (83), we conclude

$$(\mathsf{A}\gamma, \mathsf{A}\gamma') = (\gamma, \gamma')\det \mathsf{A}$$

for any $\gamma, \gamma' \in H^1$. Introducing a new variable $\gamma'' = \mathsf{A}\gamma$, we see that

$$(\mathsf{A}^{-1}\gamma'', \gamma') = \frac{1}{\det \mathsf{A}}(\gamma'', \mathsf{A}\,\gamma')$$

which is equivalent to (114).

### B.4. The general case

The torus $\Sigma = S^1 \times S^1$ is the product of two circles. We have

$$H^0(S^1, \mathbb{F}) = \mathbb{F}, \quad H^1(S^1, \mathbb{F}) = \mathbb{F}\gamma, \quad \gamma \cup \gamma = 0.$$

The relation

$$\gamma_1 \cup \gamma_2 = -\gamma_2 \cup \gamma_1$$

that we have in

$$H^\bullet(\Sigma, \mathbb{F}) = H^\bullet(S^1, \mathbb{F}) \otimes H^\bullet(S^1, \mathbb{F})$$

is an illustration of how one is supposed to put signs in the Künneth theorem (76) for odd cohomology classes.

Now we can take

$$H^\bullet((S^1)^N, \mathbb{F}) = \mathbb{F}\langle \gamma_1, \ldots, \gamma_n \rangle / (\gamma_i^2 = 0, \gamma_i\gamma_j + \gamma_j\gamma_i = 0), \tag{118}$$

where angle brackets means we don't assume that $\gamma_i$ and $\gamma_j$ commute. Indeed, they anticommute in the algebra (118). Note that the dimensions $\dim H^i((S^1)^N, \mathbb{F}) = \binom{N}{i}$ are the binomial coefficients, and hence the symmetry of the binomial coefficients may be interpreted as an instance of Poincaré duality.

When we act by $\mathsf{A}$ in the basis $\{\gamma_1, \ldots, \gamma_n\}$ of $H^1$, we get, unraveling the definitions,

$$H^N((S^1)^N, \mathbb{F}) \xrightarrow{\text{multiplication by (110)}} H^N((S^1)^N, \mathbb{F}).$$

The Poincaré duality between $H^1$ and $H^{N-1}$ gives the Cramer's formula for $\mathsf{A}^{-1}$.

## C. Tropical lines, planes, etc.
### C.1.

Consider the line

$$x_2 = ax_1 + b \subset \mathbb{C}^2, \quad a, b \neq 0, \tag{119}$$

equivalently, the graph of the function $x_1 \mapsto ax_1 + b$. What does it look like when $x_1$ and $x_2$ are exponentially large or small ? The question being a little vague, let us start by describing the set of possible values of $v_i = \ln|x_i|$ for $(x_1, x_2)$ satisfying (119).

The complex numbers $\{-x_2, ax_1, b\} \subset \mathbb{C}$ sum to zero, hence we can form a triangle in the complex planes with these vectors as sides. The triangle inequality says that a triangle with side lengths $L_1, L_2, L_3$ exists if an only if each $L_i$ is less than or equal to the sum of the other two numbers. This means

$$\begin{aligned} |x_2| &\leq |ax_1| + |b|, \\ |ax_1| &\leq |x_2| + |b|, \\ |b| &\leq |ax_1| + |x_2|, \end{aligned} \tag{120}$$

which is equivalent to

$$\ln\left(\pm|a|e^{v_1} \mp |b|\right) \leq v_2 \leq \ln\left(|a|e^{v_1} + |b|\right). \tag{121}$$

See the plot on the left in Figure (122) for $|a| = 2, |b| = 3$. Note this graph dips to $v_2 = -\infty$ precisely at the value of $v_1$ that corresponds to the unique root $x_1 = -b/a$ of $ax_1 + b = 0$.

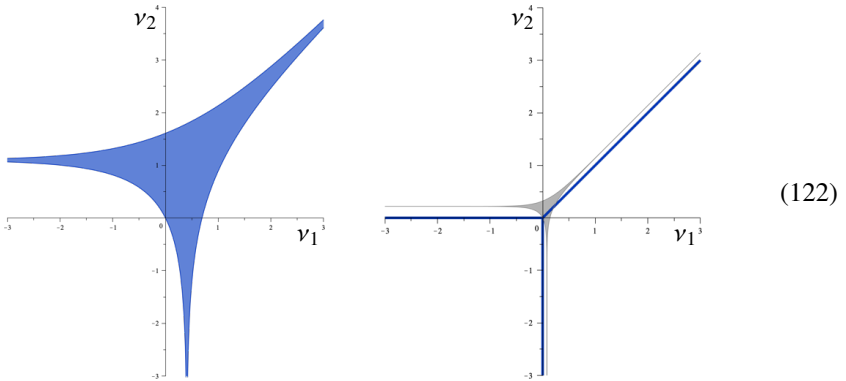

(122)

What does the plot on the left in Figure (122) look like at a very large scale ? We should rescale it by $1/T$, where $T$ some large number, and take the limit when $T \to +\infty$. For $T = 5$ we get the shape in gray on the right in Figure (122) and as $T \to \infty$ these shapes converge to the union of 3 rays plotted in blue. This union of three rays is called the *tropicalization* of the line (119).

### C.2.

Here is an alternative way to talk about $x_i$ being exponentially large or small. Above, we had a parameter $T \gg 0$, and after rescaling $v_i$ by $T$, the absolute values $|x_i|$ were of the order $e^{\nu_i T}$. To introduce objects like $e^{\nu_i T}$, where $T$ is a parameter, into the framework of linear algebra, there should be some extension of the field $\mathbb{C}$ that contains elements $e^{\xi T}$ for all $\xi \in \mathbb{R}$. Their multiplication is clear

$$e^{\xi_1 T} \cdot e^{\xi_2 T} = e^{(\xi_1 + \xi_2)T},$$

but how should we compute inverses like $(e^{\xi_1 T} + e^{\xi_2 T})^{-1}$ ?

Since we are interested in the $T \to +\infty$ limit, we should focus on the larger of the two exponents $\xi_i$ in $e^{\xi_1 T} + e^{\xi_2 T}$. Suppose $\xi_1 > \xi_2$. Then we can write

$$\frac{1}{e^{\xi_1 T} + e^{\xi_2 T}} = \frac{1}{e^{\xi_1 T}} \frac{1}{1 + e^{(\xi_2 - \xi_1)T}} = e^{-\xi_1 T} \sum_{n=0}^{\infty} (-1)^n e^{n(\xi_2 - \xi_1)T}. \quad (123)$$

The series in (123) is a geometric series which converges in the usual sense of calculus if $T$ is a positive real number. For our purposes, a much weaker notion of convergence will be sufficient.

By definition, an absolute value on an algebra $\mathbb{A}$ is a function

$$\mathbb{A} \xrightarrow{\|\cdot\|} \mathbb{R}_{\geq 0},$$

that satisfies

$$\|x\| = 0 \quad \Leftrightarrow \quad x = 0,$$
$$\|xy\| = \|x\|\|y\|, \quad (124)$$
$$\|x + y\| \leq \|x\| + \|y\|. \quad (125)$$

For example, the usual absolute value on $\mathbb{C}$ used in (120) satisfies the above conditions, and we have used the triangle inequality (125) in the derivation of (120).

Another example of an absolute value is

$$\left\|\sum c_i e^{\xi_i T}\right\|_{\sim} = e^{\xi_{\max}}, \quad \xi_{\max} = \max_{c_i \neq 0} \xi_i, \quad (126)$$

where is subscript is chosen to remind us that the absolute value (126) records the leading *asymptotics* in the $T \to +\infty$ limit. Instead of (125), this absolute value satisfies a stronger property

$$\|x + y\|_{\sim} \leq \max(\|x\|_{\sim}, \|y\|_{\sim}), \quad \text{and, moreover,} \quad (127)$$
$$\|x + y\|_{\sim} = \max(\|x\|_{\sim}, \|y\|_{\sim}) \quad \text{if } \|x\|_{\sim} \neq \|y\|_{\sim}. \quad (128)$$

Such absolute values are called *nonarchimedian*.

The series in (123) converges with respect to the absolute value (126) in the sense that

$$\left\|\sum_{n=N}^{\infty} e^{n(\xi_2 - \xi_1)}\right\|_{\sim} \to 0, \quad N \to \infty. \quad (129)$$

More generally, all series of the following form converge:

$$\mathbb{F}_{\sim} = \left\{\sum c_i e^{\xi_i T}, \text{ where } c_i \in \mathbb{C} \text{ and } \lim \xi_i = -\infty\right\}. \quad (130)$$

The reader should check that (130) is a field. The formula (126) defines an absolute value on this field. To save on notation, one can denote $t = e^{-T}$. The series (130) are then series in ascending real powers of $t$.

The unifying feature in Figure (122) is that in both cases we have the image of the line (119) under the map

$$(x_1, x_2) \mapsto (\ln \|x_1\|, \ln \|x_2\|).$$

Mathematicians call such images amoebas, because they will look a little bit like an amoeba if we replace the line by a plane curve defined by an equation of degree $\geq 3$. In other words, the tropical line is a *nonarchimedian amoeba* of a line.

### C.3.

Given an absolute value $\|\cdot\|$, we define

$$v(x) = \ln \|x\| \in \mathbb{R} \cup \{-\infty\}. \tag{131}$$

For a nonarchimedian absolute value $\|\cdot\|$, this satisfies

$$v(x) = -\infty \quad \Leftrightarrow \quad x = 0,$$
$$v(xy) = v(x) + v(y), \tag{132}$$
$$v(x+y) \in \max_Y (v(x), v(y)), \tag{133}$$

where (133) combines the two cases (127) and (128) into one formula using a multivalued function

$$\max_Y (v_1, \ldots, v_n) = \begin{cases} \max v_i, & \text{if this maximum is unique}, \\ [-\infty, \max v_i], & \text{otherwise} \end{cases} \tag{134}$$

The subscript in (134) is to remind us what the graph of this function looks like. Indeed, the graph on the right in Figure (122) is the plot of the multivalued function

$$\max_Y (0, \xi_1) = \text{possible values of } v(ax_1 + b),$$

where $v((a, b, x_1)) = (0, 0, \xi_1)$.

### C.4.

Now we are ready to generalize this discussion to a hyperplane

$$x_{n+1} = \sum_{i=1}^{n} a_i x_i + a_0, \quad (x_1, \ldots, x_{n+1}) \in \mathbb{C}^{n+1}, \tag{135}$$

where all coefficients $a_i$ are nonzero complex numbers. All arguments above generalize verbatim and give

$$\xi_{n+1} = \max_Y (0, \xi_1, \ldots, \xi_n) \tag{136}$$

as the tropicalization of (135). This is what the plot of this function looks like for $n = 2$. This is a tropical hyperplane in 3-space.

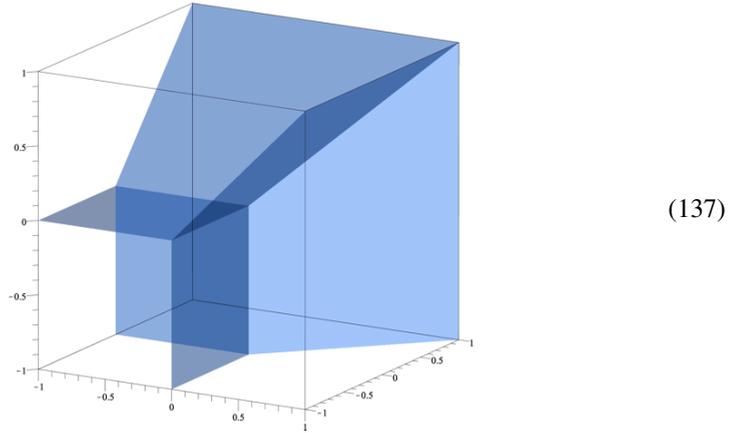

(137)

Note that we have

$$\max\nolimits_\Upsilon (0, \xi_1, \xi_2) = -\infty = \nu(0), \qquad (138)$$

precisely for $(\xi_1, \xi_2)$ forming a tropical line. This is a tropical analog of the obvious fact the the intersection of $x_{n+1} = 0$ with another hyperplane is a hyperplane in $\mathbb{C}^n$.

### C.5.

In place of a linear polynomial in (135), we could have taken an arbitrary polynomial

$$x_{n+1} = P(x_1, \ldots, x_n) = \sum_{\boldsymbol{\beta} \in \mathbb{Z}^n} a_{\boldsymbol{\beta}} x^{\boldsymbol{\beta}}, \quad a_{\boldsymbol{\beta}} \in \mathbb{F}_\sim, \qquad (139)$$

where $\boldsymbol{\beta} = [\beta_1, \ldots, \beta_n] \in \mathbb{Z}^n$ and

$$x^{\boldsymbol{\beta}} = \prod_{i=1}^n x_i^{\beta_i}.$$

In (139) we assume that only finitely many coefficients $a_{\boldsymbol{\beta}}$ are nonzero. Note that here we allow $a_{\boldsymbol{\beta}}$ to be any elements of $F_\sim$. In other words, we allow the coefficients to be exponentially large or small.

Arguing as above, we get

$$\xi_{n+1} = \max\nolimits_\Upsilon \left\{ \langle \boldsymbol{\beta}, \boldsymbol{\xi} \rangle + \nu(a_{\boldsymbol{\beta}}) \right\}_{a_{\boldsymbol{\beta}} \neq 0} \qquad (140)$$

as the tropicalization of (139). Here $\boldsymbol{\xi} = [\xi_1, \ldots, \xi_n]^{\mathrm{T}}$ and the angle brackets were defined in (104).

The left-hand side in (140) is called a tropical polynomial[16]. The set (140) in $\mathbb{R}^{n+1}$ is the graph of this polynomial. The intersection of the graph with $\xi_{n+1} = -\infty$ is the tropicalization of the hypersurface $P(x_1, \ldots, x_n) = 0$. Here is an example of a graph of a tropical

---

[16] Real numbers $\mathbb{R}$ with operations $\{\max, +\}$ form a semiring called tropical semiring

polynomial of degree 3 in two variables:

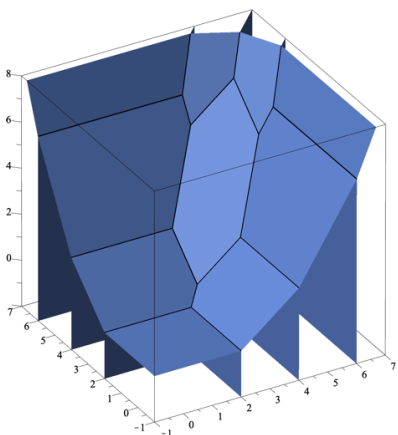

(141)

The facets in (141) have slopes $\boldsymbol{\beta} \in \{[0,0], [1,0], \ldots, [3,0], [2,1], [1,2], [0,3]\}$ and they move up and down as $v(a_\beta)$ change, leading to changes in the combinatorics. The reader should experiment to get some feeling for how this works. Graphs of linear polynomials from (122) and (137) change only by an overall translation if we take the coefficients from $\mathbb{F}_\sim$ instead of $\mathbb{C}$.

### C.6.

Tropical varieties of codimension more than 1 are certain piecewise linear objects that are defined axiomatically. See [27, 33]. In particular, they *do not* need to be nonarchimedian amoebas of any an algebraic variety over $\mathbb{F}_\sim$ or some other field with a nonarchimedian norm.

In fact, this already happens for linear spaces. There is a tropical linear space for any matroid including the nonrealizable ones ! To see how it works, let us go back to the settings of Section 5.3.

Let us first consider the case when the $Y \subset \mathbb{F}^n$ is a linear subspace, where $\mathbb{F}$ is a field. Given a subset $S \subset \{1, \ldots, n\}$, let $\mathbb{F}^S$ denote the quotient space of $\mathbb{F}^n$ with the coordinates $y_i$, where $i \in S$. We have

$$|S| = \text{rank } S \quad \Leftrightarrow \quad \begin{array}{c} \text{the map } Y \to \mathbb{F}^S \\ \text{is surjective} \end{array}. \tag{142}$$

For any matroid, the subsets $S$ such that $|S| = \text{rank } S$ are called *independent*. Minimal dependent subsets are called *circuits*. One can reconstruct the matroid completely from the knowledge of either the independent sets or the circuits.

For a circuit $S$, the image of $Y \to \mathbb{F}^S$ is a hyperplane with an equation that involves all variables $y_i$ with $i \in S$. We already know how to tropicalize it. One should take the tropical hyperplane

$$\max_Y(\{y_i\}_{i \in \text{circuit } S}) = -\infty. \tag{143}$$

Very remarkably, for an *arbitrary* matroid $M$, the equations (143), taken for all circuits $S$ of $M$, describe a tropical linear space. This linear space provides a tropical realization of the matroid.

## References

[1] Omid Amini and Matthieu Piquerez, *Hodge theory for tropical varieties*. arxiv.org/abs/2007.07826. ↑13

[2] ———, *Tropical Clemens-Schmid sequence and existence of tropical cycles with a given cohomology class*. arxiv.org/abs/2012.13142. ↑13

[3] Federico Ardila and Adam Boocher, *The closure of a linear space in a product of lines*, J. Algebraic Combin. **43** (2016), no. 1, 199–235. ↑13

[4] Sheldon Axler, *Linear algebra done right*, 3rd ed., Undergraduate Texts in Mathematics, Springer, Cham, 2015. ↑28

[5] Matthew Baker, *Hodge theory in combinatorics*, Bull. Amer. Math. Soc. (N.S.) **55** (2018), no. 1, 57–80. ↑28

[6] Matthew Baker and Nathan Bowler, *Matroids over partial hyperstructures*, Adv. Math. **343** (2019), 821–863. ↑13

[7] Matthew Baker and Oliver Lorscheid, *The moduli space of matroids*, Adv. Math. **390** (2021), Paper No. 107883, 118. ↑13

[8] ———, *Foundations of matroids I: Matroids without large uniform minors*. arxiv.org/abs/2008.00014. ↑13

[9] Anders Björner and Torsten Ekedahl, *On the shape of Bruhat intervals*, Ann. of Math. (2) **170** (2009), no. 2, 799–817. ↑5

[10] Tom Braden, June Huh, Jacob Matherne, Nicholas Proudfoot, and Botong Wang, *Singular Hodge theory for combinatorial geometries*. arxiv.org/abs/2010.06088. ↑8, 9, 18, 19, 20, 21, 28

[11] Francesco Brenti, *Unimodal, log-concave and Pólya frequency sequences in y combinatorics*, Mem. Amer. Math. Soc. **81** (1989), no. 413. ↑17

[12] N. G. de Bruijn and P. Erdös, *On a combinatorial problem*, Nederl. Akad. Wetensch., Proc. **51** (1948), 1277–1279 = Indagationes Math. 10, 421–423 (1948). ↑3

[13] Jordana Cepelewicz, Quanta Magazine. to appear. ↑28

[14] R. P. Dilworth and Curtis Greene, *A counterexample to the generalization of Sperner's theorem*, J. Combinatorial Theory Ser. A **10** (1971), 18–21. ↑7

[15] Thomas A. Dowling and Richard M. Wilson, *The slimmest geometric lattices*, Trans. Amer. Math. Soc. **196** (1974), 203–215. ↑5

[16] ———, *Whitney number inequalities for geometric lattices*, Proc. Amer. Math. Soc. **47** (1975), 504–512. ↑5

[17] Ben Elias and Geordie Williamson, *The Hodge theory of Soergel bimodules*, Ann. of Math. (2) **180** (2014), no. 3, 1089–1136. ↑5, 21

[18] Alex Fink, *Tropical cycles and Chow polytopes*, Beitr. Algebra Geom. **54** (2013), no. 1, 13–40. ↑13

[19] Curtis Greene, *A rank inequality for finite geometric lattices*, J. Combinatorial Theory **9** (1970), 357–364. ↑7

[20] Kevin Hartnett, *A Path Less Taken to the Peak of the Math World*, Quanta Magazine ( June 27, 2017). www.quantamagazine.org/a-path-less-taken-to-the-peak-of-the-math-world-20170627/. ↑28

[21] Allen Hatcher, *Algebraic topology*, Cambridge University Press, Cambridge, 2002. ↑28

[22] June Huh and Botong Wang, *Enumeration of points, lines, planes, etc*, Acta Math. **218** (2017), no. 2, 297–317. ↑6, 28

[23] June Huh, *Combinatorial applications of the Hodge-Riemann relations*, Proceedings of the International Congress of Mathematicians—Rio de Janeiro 2018. Vol. IV. Invited lectures, World Sci. Publ., Hackensack, NJ, 2018, pp. 3093–3111. ↑28

[24] ———, *Combinatorics and Hodge theory*. Proceedings of the International Congress of Mathematicians 2022. ↑5, 28

[25] A. W. Ingleton, *Representation of matroids*, Combinatorial Mathematics and its Applications (Proc. Conf., Oxford, 1969), Academic Press, London, 1971, pp. 149–167. ↑13

[26] A. W. Ingleton and R. A. Main, *Non-algebraic matroids exist*, Bull. London Math. Soc. **7** (1975), 144–146. ↑13

[27] Ilia Itenberg, Grigory Mikhalkin, and Eugenii Shustin, *Tropical algebraic geometry*, 2nd ed., Oberwolfach Seminars, vol. 35, Birkhäuser Verlag, Basel, 2009. ↑13, 41

[28] Gil Kalai, *The Work of June Huh*. Proceedings of the International Congress of Mathematicians 2022. ↑28

[29] Kalle Karu, *Hard Lefschetz theorem for nonrational polytopes*, Invent. Math. **157** (2004), no. 2, 419–447. ↑21

[30] Eric Katz, *What is . . . tropical geometry?*, Notices Amer. Math. Soc. **64** (2017), no. 4, 380–382. ↑13

[31] David Kazhdan and George Lusztig, *Schubert varieties and Poincaré duality*, Geometry of the Laplace operator (Proc. Sympos. Pure Math., Univ. Hawaii, Honolulu, Hawaii, 1979), Proc. Sympos. Pure Math., XXXVI, Amer. Math. Soc., Providence, R.I., 1980, pp. 185–203. ↑21

[32] Bernt Lindström, *The non-Pappus matroid is algebraic*, Ars Combin. **16** (1983), no. B, 95–96. ↑13

[33] Diane Maclagan and Bernd Sturmfels, *Introduction to tropical geometry*, Graduate Studies in Mathematics, vol. 161, American Mathematical Society, Providence, RI, 2015. ↑13, 41

[34] George Melvin and William Slofstra, *Soergel bimodules and the shape of Bruhat intervals*, preprint, 2020. pages 5

[35] Th. Motzkin, *Beiträge zur Theorie der linearen Ungleichungen*, 1936. Dissertation Thesis, University of Basel, Jerusalem. ↑3

[36] ———, *The lines and planes connecting the points of a finite set*, Trans. Amer. Math. Soc. **70** (1951), 451–464. ↑3

[37] Peter Nelson, *Almost all matroids are nonrepresentable*, Bull. Lond. Math. Soc. **50** (2018), no. 2, 245–248, DOI 10.1112/blms.12141. MR3830117 ↑12

[38] James G. Oxley, *Matroid theory*, Oxford Science Publications, The Clarendon Press, Oxford University Press, New York, 1992. ↑28

[39] Gian-Carlo Rota, *Combinatorial theory, old and new*, Actes du Congrès International des Mathématiciens (Nice, 1970), Gauthier-Villars, Paris, 1971, pp. 229–233. ↑5

[40] Igor R. Shafarevich, *Basic algebraic geometry. 1*, Translated from the 2007 third Russian edition, Springer, Heidelberg, 2013. ↑28

[41] Richard P. Stanley, *The number of faces of a simplicial convex polytope*, Adv. in Math. **35** (1980), no. 3, 236–238. ↑21

[42] ———, *Combinatorial applications of the hard Lefschetz theorem*, Proceedings of the International Congress of Mathematicians, Vol. 1, 2 (Warsaw, 1983), PWN, Warsaw, 1984, pp. 447–453. ↑28

[43] ———, *Algebraic combinatorics*, Undergraduate Texts in Mathematics, Springer, Cham, 2018. Walks, trees, tableaux, and more; Second edition of [ MR3097651]. ↑28

[44] Bernd Sturmfels, *Solving systems of polynomial equations*, CBMS Regional Conference Series in Mathematics, vol. 97, Published for the Conference Board of the Mathematical Sciences, Washington, DC; by the American Mathematical Society, Providence, RI, 2002. ↑13

### Andrei Okounkov

Andrei Okounkov, Department of Mathematics, University of California, Berkeley, 970 Evans Hall Berkeley, CA 94720–3840, okounkov@math.columbia.edu